\documentclass[12pt,twoside]{article}
\usepackage{amsmath,amsthm,amssymb}
\usepackage{amsfonts}

\newenvironment{Proof}{\textbf{Proof.}}{$\qquad \blacksquare$\par}
\newenvironment{Proof of}[1]{\textbf{Proof of the #1.}}{$\qquad \blacksquare$\par}
\newenvironment{Item}[1]{\par #1 }{\par}

\newcommand{\B}{\mathcal B}
\newcommand{\CC}{\mathcal C}

\newcommand{\x}{\widetilde x}
\newcommand{\y}{\widetilde y}
\newcommand{\X}{\widetilde X}

\newcommand{\U}{\widetilde U}
\newcommand{\V}{\widetilde V}
\newcommand{\F}{\widetilde F}
\newcommand{\TDelta}{\widetilde \Delta}
\newcommand{\Sol}{\mathcal S}
\newcommand{\Tdelta}{\widetilde \delta}
\newcommand{\tal}{\widetilde \alpha}
\newcommand{\al}{\alpha}

\newcommand{\A}{\mathcal A}

\newcommand{\C}{\mathbb C}
\newcommand{\R}{\mathbb R}
\newcommand{\Z}{\mathbb Z}
\newcommand{\N}{\mathbb N}

\newcommand{\TT}{\mathcal T}

\newtheorem{thm}{Theorem}[section]
\newtheorem{lem}[thm]{Lemma}
\newtheorem{prop}[thm]{Proposition}
\newtheorem{corl}[thm]{Corollary}

\theoremstyle{definition}
\newtheorem{defn}[thm]{Definition}
\newtheorem{exm}[thm]{Example}
\newtheorem{rem}[thm]{Remark}

\numberwithin{equation}{section}

\begin{document}
\thispagestyle{empty}
\begin{flushright}
\scriptsize{IM UwB$\qquad$\\
July  2004$\qquad$}
\end{flushright}
\begin{center}
{\bf \LARGE Covariance algebra of a partial dynamical system}\\
\bigskip
{\large \textsc{B. K.  Kwa\'sniewski} }\\
\bigskip
  Institute of Mathematics,  University  in Bialystok,\\
 ul. Akademicka 2, PL-15-424 Bialystok, Poland\\
 e-mail: bartoszk@math.uwb.edu.pl
\end{center}
\begin{abstract} 
\noindent A pair $(X,\al)$ is a partial dynamical system if $X$
is a compact topological space and $\al:\Delta \rightarrow X$ is a
continuous mapping such that $\Delta$ is open. Additionally we
assume here that $\Delta$ is closed and $\al(\Delta)$ is open. Such systems arise naturally while dealing with commutative
$C^*$-dynamical systems.
\\
 In this paper  we construct and
investigate a universal $C^*$-algebra $C^*(X,\al)$ which agrees
with the
  partial crossed product  \cite{exel} -  in the case $\al$ is injective,
   and   with the  crossed product by a monomorphism \cite{Murphy} - in the case $\al$ is onto.
\\
   The main method here is to use the description  of maximal ideal
space of a coefficient algebra, cf. \cite{maxid},
\cite{lebiedodzij}, in order  to construct a larger system
$(\X,\tal)$ where $\tal$ is a partial homeomorphism. Hence one may
apply in extenso the  partial crossed product theory \cite{exel},
\cite{exel_laca_quigg}.   In particular, one generalizes the
notions of topological freedom and invariance of a set, which are
being used afterwards to obtain The Isomorphism Theorem and the
complete description of ideals of $C^*(X,\al)$.
\end{abstract}


\tableofcontents

% D O U B L E   S P A C I N G
% To TURN ON DOUBLE SPACING remove sign "%" from the next line
%\addtolength{\baselineskip}{13pt.}

\section*{Introduction}
\addcontentsline{toc}{section}{Introduction}
In  quantum theory  the term covariance algebra (crossed product) means  an algebra generated by an algebra of observables and by operators which determine the time evolution of a quantum system (a $C^*$-dynamical system), thereby  the covariance algebra is an object which carries all the information about the quantum system, see 
\cite{Brat_Robin}, \cite{lebiedodzij} (and the sources cited there) for this and other connections with mathematical physics. In pure mathematics $C^*$-algebras associated to $C^*$-dynamical systems  proved to be useful in different fields:  classification of operator algebras \cite{Pedersen}, \cite{Brat_Robin}, \cite{Odonovan}; K-theory for $C^*$-algebras \cite{blackadar}, \cite{exel}, \cite{mcclanachan}; functional and functional differential equations \cite{Anton}, \cite{Anton_Lebed}; or even in number theory \cite{lr}.  This multiplicity of applications and the complexity of the matter attracted many authors and caused  an abundence of various  approaches  \cite{Pedersen}, \cite{exel}, \cite{mcclanachan},  \cite{Stacey}, \cite{Adji_Laca_Nilsen_Raeburn}, \cite{Murphy}, \cite{exel2}, \cite{exel_royer}. In the present article we propose another approach which on one hand may  seem  to embrace  a very special case but on the other hand: 
\begin{itemize}
\item[1)] unlike the other authors investigating crossed products (see discussions below) we do not require here any kind (substitute) of reversibility of the given system,
\item[2)] we  obtain a rather thorough description of the associated covariance algebra, and also strong tools to study it, 
\item[3)] we find points of contact of  different approaches and thereby  clarify the relations between them. 
\end{itemize}
\medskip 

 In order to give the motivation of the  construction of the crossed product developed in the paper  
we would like to present  a simple example.
\ \\

\textbf{Example.}
Consider the Hilbert space $H=L^2_\mu({\mathbb R})$ where $\mu$ is the Lebesgue measure. Let
$\A\subset L(H)$ be the $C^*$-algebra of operators of
multiplication by continuous bounded functions on $\mathbb R$ that
are constant on ${\mathbb R}_- = \{ x: x\le 0 \}$. Set the unitary
operator $U\in L(H)$ by the formula
$$
(Uf)(x)= f (x-1), \qquad f (\,\cdot\, )\in H .
$$
Routine verification shows that the mapping
$$
\A\ni a \mapsto UaU^*
$$
is an  endomorphism of $\A$ of the form
\begin{equation}\label{U3}
UaU^* (x) = a(x-1), \qquad a(\,\cdot\, ) \in \A ,
\end{equation}
and
\begin{equation}\label{U4}
U^*aU (x) = a(x+1), \qquad a(\,\cdot\, ) \in \A .
\end{equation}
Clearly the mapping $\A\ni a \mapsto U^*aU$ {\em does not} preserve $\A$.

Let $C^* (\A,U)$ be the $C^*$-algebra generated by $\A$ and $U$. It is
easy to show that
$$
C^* (\A,U) = C^* (\B,U)
$$
where $\B\subset L(H)$ is the $C^*$-algebra of operators of
multiplication by continuous bounded functions on $\R$ that
have limits at $-\infty$.

In addition we have that
$$
 U\B U^*\subset \B \qquad \text{and}\qquad
U^* \B U\subset \B
$$
and the corresponding actions $\delta (\,\cdot\, )= U(\,\cdot\, )U^*$
and $\delta_* (\,\cdot\, )= U^* (\,\cdot\, )U$ on $\B$ are given
by formulae \eqref{U3} and \eqref{U4}.

Moreover 
$$
C^* (\A,U) = C^* (\B,U)\cong \B\times_{\delta} \Z.
$$
Where in the right hand side stands the standard crossed product of $\B$ by the {\em automorphism} $\delta$.

This example shows a natural situation when the crossed product $\B\times_{\delta}\Z$ is 'invisible' at the beginig  (on the initial algebra $\A$, $\delta$ acts as an endomorphism and $\delta_*$ even does not preserve $\A$) but after implementing a natural extension of $\A$ up to $\B$,  $\delta$ becomes an {\em automorphism} and leads to the 
crossed product structure. The aim of the paper is to investigate the general constructions of this type.  
We will also find out that the arising constructions are rather natural and one can come across 
them 'almost anywhere', in particular they present the passage from the irreversible topological Markov chains to the reverible ones (see Propostition \ref{topologicalmarkovchain}) and the maximal ideal spaces of the algebras of $\B$ type possess the solenoid  structure (see  Examples \ref{solenoid}, \ref{solenoid2}).
\ \\

We   deal here with $C^*$-dynamical systems where dynamics is implemented
by a single endomorphism, hence  a $C^*$-\emph{dynamical system}
is identified with a pair $(\A,\delta)$ where $\A$ is a unital
$C^*$-algebra and $\delta:\A\rightarrow\A$ is a $^*$-endomorphism. Additionally we  assume that
 $\A$ is  commutative.  Hence, in fact, we deal  with topological dynamical systems. Indeed, using the Gelfand transform one can identify  $\A$ with the algebra $C(X)$
   of continuous functions  on  the maximal ideal
space $X$ of $\A$, and within this identification the endomorphism
$\delta$ generates (see, for example, \cite{maxid}) a
continuous \emph{partial mapping}   $\alpha:\Delta \rightarrow X$
where $\Delta\subset X$ is closed and open (briefly clopen), and
the following formula holds
 \begin{equation} \label{e3.0}
\delta (a)(x)=\left\{ \begin{array}{ll} a(\alpha(x))& ,\ \ x\in
\Delta\\
0 & ,\ \ x\notin \Delta \end{array}\right. ,\qquad a\in C(X).
\end{equation}
Therefore  we have  one-to-one correspondence between the
commutative unital $C^*$-dynamical systems $(\A,\delta)$ and pairs
$(X,\al)$, where $X$ is compact and $\al$ is a partial continuous
mapping which domain is clopen. We shall call  $(X,\al)$ a
\emph{partial dynamical system}.

Usually  covariance
algebra is  another name for the crossed product which in turn is defined in various ways \cite{Pedersen}, \cite{Odonovan}, \cite{exel}, \cite{Murphy}, though it seems
more appropriate to define such objects as  $C^*$-algebras with a universal
property with  respect to covariant representations \cite{Stacey}, \cite{Adji_Laca_Nilsen_Raeburn}. In
the literature, cf. \cite{Pedersen}, \cite{exel},
\cite{mcclanachan}, \cite{Adji_Laca_Nilsen_Raeburn},  a
\emph{covariant representation} of $(\A,\delta)$   is  meant to be
a triple $(\pi,U,H)$   where  $H$ is a Hilbert space,
$\pi:\A\rightarrow L(H)$ is a representation of $\A$ by bounded
operators  on $H$, and $U\in L(H)$ is such that
$$
\pi(\delta(a))=U \pi(a) U^*,\qquad
\mathrm{for\,\,all}\,\,\,\,a\in\A,
$$
  plus eventually some  other  conditions   imposed on $U$ and $\pi$. If
$\pi$ is faithful we shall  call $(\pi,U,H)$ a \emph{covariant
faithful representation}. 
 The covariant representations of a
$C^*$-dynamical system give rise to a category
$\mathrm{Cov}(\A,\delta)$ where objects  are  the $C^*$-algebras
$C^*(\pi(\A),U)$, generated by $\pi(\A)$ and $U$,  while morphisms
are the usual $^*$-morphisms $\phi:C^*(\pi(\A),U)\rightarrow
C^*(\pi'(\A),U')$ such that
$$
\phi(\pi(a))=\pi'(a),\,\,\mathrm{for}\,\,
a\in\A,\qquad\mathrm{and}\qquad\phi(U)=U'
$$
 (here $(\pi,U,H)$ and $(\pi',U',H')$ denote  covariant representations of $(\A,\delta)$).
 In many cases the main interest is concentrated on  the subcategory $\mathrm{CovFaith}(\A,\delta)$
  of $\mathrm{Cov}(\A,\delta)$ for  which objects are  algebras $C^*(\pi(\A),U)$ where now $
  \pi$ is  faithful.
 The fundamental  problem then is  to describe a universal object in $\mathrm{Cov}(\A,\delta)$,
 or in $\mathrm{CovFaith}(\A,\delta)$,  in terms of the $C^*$-dynamical system $(\A,\delta)$.
 If such  an object exists then it is unique  up to isomorphism, and it shall be  called  a \emph{covariance algebra}.
 \par
 It  is well known  \cite{Pedersen} that,  in case $\delta$ is an
automorphism, the classic crossed product $\A\rtimes_\delta \Z$
is the covariance algebra of the $C^*$-dynamical system
$(\A,\delta)$. Being motivated by the paper \cite{cuntz}, in which
J. Cuntz discussed a concept of the crossed product by an
endomorphism which is not automorphism, many authors proposed theories of generalized
crossed products with some kind of universality (see
\cite{Paschke}, \cite{exel}, \cite{mcclanachan}, \cite{Murphy},
\cite{exel2}). For instance,   G. Murphy
in \cite{Murphy} has proved  that a corner of the crossed product
of a certain direct limit is a covariance algebra of a system
$(\A,\delta)$ where $\delta$ is a monomorphism (in fact he has
proved far more general result, see \cite[Theorem 2.3]{Murphy}).
R. Exel in \cite{exel} introduced  a partial crossed product which
can be applied also in the case $\delta$ is not injective, though
generating a partial automorphism (see also
\cite{exel_laca_quigg}, \cite{mcclanachan}). Nevertheless, in
general the inter-relationship between the $C^*$-dynamical system
and its covariance algebra is still not totally-established.
\par
In the  approach developed in this paper we explore the leading
concept of  the coefficient  algebra, introduced in
\cite{lebiedodzij}.  The  elements of this algebra    play the
role of Fourier's coefficients in the covariance algebra,  hence
the name.
     The authors of \cite{lebiedodzij} studied the $C^*$-algebra $C^*(\A,U)$ generated by
     a $^*$-algebra $\A\subset L(H)$ and a partial isometry $U\in L(H)$.
They have defined $\A$ (in a slightly different yet equivalent
form) to be a \emph{coefficient algebra} of $C^*(\A,U)$ whenever
$\A$ possess the following three properties
\begin{equation}
\label{b1}
 U^*U\in \A',
 \qquad
 U\A U^*\subset A,
 \qquad
 U^*\A U\subset \A,
\end{equation}
where $\A'$ denotes  the commutant of $\A$.  Let us indicate that
this concept appears, in more or less explicit form, in all
aforesaid articles: If $U$ is unitary then (\ref{b1}) holds iff
$\delta(\cdot)=U(\cdot)U^*$ is an automorphism of $\A$, and thus
in this case  $\A$ can be regarded as a coefficient algebra of the
crossed product $\A\rtimes_\delta \Z$. For example in the paper
\cite{cuntz}, the UHF algebra $\mathcal{F}_n$  is a coefficient
algebra of the Cuntz algebra $\mathcal{O}_n$. Also the algebra $A$
considered by Paschke in \cite{Paschke} is a coefficient algebra
of the $C^*$-algebra $C^*(A,S)$ generated by $A$ and the isometry
$S$. The algebra $C^{\tilde{\al}}$ defined in \cite{Murphy} as the
fixed point algebra for dual action  can be thought of as a
generalized coefficient algebra of the crossed product
$C^*(A,M,\al)$ of $A$ by the semigroup $M$ of injective
endomorphisms.
\\
Thanks to \cite{maxid}, the main tool we are given is  the
description of
   maximal ideal spaces of certain coefficient algebras. More precisely, for any partial isometry $U$
   and unital commutative $C^*$-algebra $\A$ such that $U^*U\in \A'$ and $U\A U^*\subset \A$ we
   infer
   that  $(\A,\delta)$ is a $C^*$-dynamical system, where  $\delta(\cdot)=U(\cdot)U^*$.
    However $\A$ does not need to fulfill the third property from (\ref{b1}).
     The solution then is to  pass to a bigger $C^*$-algebra $\B$
     generated by $\{\A, U^*\A U,U^{2*}\A U^2,...\}$. Then $(\B,\delta)$ is a
     $C^*$-dynamical system and $\B$ is a coefficient algebra of
     $C^*(\B,U)=C^*(\A,U)$, see \cite{lebiedodzij}. In this case
the authors of \cite{maxid} managed to 'estimate' the maximal ideal
space $M(\B)$ of $\B$ in terms of
$(\A,\delta)$, or better to say, in terms of the generated partial
dynamical system  $(X,\al)$. Fortunately, the full description of
$M(\B)$ is obtained \cite[3.4]{maxid} by a slight
strengthening of assumptions -  namely by assuming that the projection
$U^*U$ belongs not only to commutant $\A'$ but to $\A$ itself.
  The  partial dynamical system $(M(\B),\tal)$,  corresponding to $(\B,\delta)$,
   is thus completely determined by $(X,\al)$. Two important facts are to be noticed: $\tal$ is a partial homeomorphism,
    and  $U^*U\in\A$ implies that the image $\al(\Delta)$ of the partial mapping $\al$ is open, see Section 1 for details.
    \par
In Section 2, to   an arbitrary   partial dynamical
system $(X,\al)$ such that $\al(\Delta)$ is open   we associate
 another partial dynamical system $(\X,\tal)$ such that:
 \\[6pt]
 1) $\tal$ is a partial homeomorphism,
 \\[6pt]
  2)
  there exist  a continuous surjection $\Phi:\X\rightarrow X$ such that the equality $\Phi\circ \tal=\al\circ\Phi$
   holds wherever it makes sense (see  diagram
  (\ref{phidiagram})),
  \\[6pt]
  3) if $\al$ is injective then $\Phi$ becomes a homeomorphism, that
is $(\X,\tal) \cong (X,\al)$.
\\[6pt]
This authorizes  us to call
$(\X,\tal)$ the reversible extension of $(X,\al)$. In the case
$\al$ is onto, $\X$ is a projective limit (see Proposition
\ref{projectivelimit}) and thus $(\X,\tal)$ generalizes the known
construction.
\par
In Section 3  we find out that all the objects of
$\mathrm{CovFaith}(\A,\delta)$ have
 the same (up to  isomorphism) coefficient $C^*$-algebra  whose  maximal ideal space is
 $\X$.
We denote this $C^*$-algebra by $\B$. Then we
construct a coefficient $^*$-algebra $\B_0$ (the closure of $\B_0$ is
$\B$)
   with the help of which we express the interrelations between the covariant representations of $(\A,\delta)$
  and $(\B,\Tdelta)$ where $\Tdelta$ is an endomorphism associated to the partial homeomorphism
  $\tal$.
 In particular  we show  that, if $\delta$ is injective, or equivalently $\al$ is onto,
   then we have  a natural one-to-one correspondence between aforementioned representations.
    In general this correspondence is maintained only if we constrain ourselves to covariant faithful representations.
\par
In Section 4 we define $C^*(\A,\delta)=C^*(X,\al )$ to be the
partial crossed product of $\B=C(\X)$ by a partial automorphism
generated by the partial homeomorphism $\tal$. We show that
$C^*(\A,\delta)$ is the universal object in
$\mathrm{CovFaith}(\A,\delta)$, and in case $\delta$ is injective,
it is also universal when considered as an object of
$\mathrm{Cov}(\A,\delta)$. Therefore we call it a covariance
algebra.
\par
Section 5 is devoted to two important notions in $C^*$-dynamical
systems theory: topological freedom and invariant sets.
Classically, these notions  were related only to homeomorphisms,
but recently they have been adopted (generalized), by authors of
\cite{exel_laca_quigg}, to work with partial homeomorphisms, see
also \cite{lebied}. Inspired by this line of development we
present here the definitions of  topological freedom and
invariance under  a partial mapping which include also noninjective
partial  mappings.  Let us mention that, for instance, in
\cite{exel_vershik}  appears  also the  definition of
topologically free irreversible dynamical system, but the authors
of \cite{exel_vershik} attach to dynamical systems
\textit{different} $C^*$-algebras than we do, hence they are in
need of a different definition. We show that there exists  a
natural bijection between closed $\al$-invariant subsets of $X$
and closed $\tal$-invariant subsets of $\X$ and that the partial
dynamical system $(X,\al)$ is topologically free if and only if
its reversible extension $(\X,\tal)$ is topologically free.
  \par
Section  6 contains two important results. Namely, we
establish a one-to-one correspondence between the ideals in
$C^*(X,\al)$ and closed invariant subsets of $X$ generalizing
Theorem 3.5 from \cite{exel_laca_quigg}. Then we present also a
version of the Isomorphism Theorem, cf.
\cite[Theorem 7.1]{Anton}, \cite[Chapter 2]{Anton_Lebed}, \cite[Theorem 2.6]{exel_laca_quigg},
\cite[Theorem 2.13]{lebiedodzij}, which
says that all objects of $\mathrm{CovFaith}(\A,\delta)$ are
isomorphic to $C^*(\A,\delta)$ whenever the corresponding system
$(X,\al)$ is topologically free.
 \section{Preliminaries. Maximal ideal space of  a  coefficient $C^*$-algebra}
 We start this section with  fixing of some notation. Afterwards,   we    present and discuss briefly the
results of \cite{maxid} in order to present our
methods and motivations. We finish this section with Theorem
\ref{tu3} to be used extensively in the sequel.
\par
 Throughout this article $\A$ denotes a
commutative unital $C^*$-algebra,   $X$ denotes its maximal ideal
space (i.e. a compact topological space), $\delta$ is an
endomorphism of $\A$, while  $\al$ stands for a continuous partial
mapping $\alpha:\Delta \rightarrow X$ where $\Delta\subset X$ is
clopen and the  formula (\ref{e3.0}) holds. We adhere to the
convention that $\N=0,1,2,...$, and when dealing with partial
mappings we follow the notation of \cite{maxid}, i.e.: for $n>0$,
we denote the domain of $\al^n$  by $\Delta_n=\al^{-n}(X)$ and its
image by $\Delta_{-n}=\al^{n}(\Delta_n)$; for $n=0$, we set
$\Delta_0=X$ and thus, for $n,m\in \N$,  we have
\begin{equation}
 \label{b-2}
 \al^n:\Delta_n\rightarrow  \Delta_{-n},
\end{equation}
 \begin{equation}\label{b-3} \al^n (\al^m(x)) = \al^{n+m}(x),\qquad x\in
\Delta_{n+m}.
\end{equation}
We recall  that in terms of the multiplicative  functionals of
$\A$, $\al$ is given by
 \begin{equation} \label{e1.0}
 x\in \Delta_1 \Longleftrightarrow \, x(\delta(1))=1, \end{equation}
 \begin{equation}\label{e2.0}
 \al(x)=x\circ \delta,\qquad x \in\Delta_1.
\end{equation}
 For the purpose of the present section we  fix (only in this section) a
  faithful representation of $
\A$, i.e. we assume  that $\A$ is a  $C^*$-subalgebra of $L(H)$
   where $L(H)$ is an  algebra of bounded linear operators on a Hilbert space $H$.
    Additionally we assume that endomorphism  $\delta$ is given by
the formula
$$
\delta(a)= UaU^*,\qquad a\in \A,
$$
for some $U\in L(H)$ and so   $U$ is a  partial isometry (note
that there exists  a  correspondence between properties of $U$ and
the partial mapping $\al$, cf.  \cite[2.4]{maxid}). In that case,
as it makes sense, we will  consider $\delta$ also as a mapping
on $L(H)$. There is a point in studying  together with
$\delta(\cdot)=U(\cdot)U^*$ one more  mapping
$$
\delta_*(b)=U^*bU, \qquad b\in L(H),
$$
which in general  maps $a\in \A$ onto an element outside  the
algebra $\A$ and hence, even if we assume that $U^*U\in
 \A'$,  we need
to pass to a bigger algebra  in order to
obtain an algebra satisfying (\ref{b1}) .
\begin{prop}\emph{\cite[Proposition 4.1]{lebiedodzij}}\label{prop1}
If $\delta(\cdot)=U(\cdot)U^*$ is an endomorphism of $\A$,
$U^*U\in
 \A'$ and $\B=C^*(\bigcup_{n=0}^\infty U^{*n}\A U^n)$
 is a $C^*$-algebra generated by
$\bigcup_{n=0}^\infty U^{*n} \A U^n$, then $\B$ is
commutative and both the mappings
$\delta:\B\rightarrow \B$ and
$\delta_*:\B\rightarrow \B$ are
endomorphisms.
\end{prop}
\noindent The elements of the algebra $\B$ play
the role of coefficients in a $C^*$-algebra $C^*(\A,U)$ generated
by $\A$ and $U$,  \cite[2.3]{lebiedodzij}. Hence the authors of
\cite{lebiedodzij} call $\B$  a \emph{coefficient
algebra}.
 It is of primary importance  that $\B$ is commutative and that we
have a description  of its  maximal ideal space,
denoted  here  by $M(\B)$,  in terms of the
maximal ideals in $\A$, see  \cite{maxid}.  Let us  recall it.
\\
\noindent With every $\widetilde{x}\in M(\B)$ we associate  a sequence of functionals
$\xi^n_{\widetilde{x}}:{\cal A}\rightarrow \C\,$, $n\in\N$,
defined by the condition
\begin{equation} \label{xi2}
\xi^n_{\widetilde{x}}(a)=\delta_*^n(a)(\widetilde{x}),\quad a\in
\A.
\end{equation}
The sequence $\xi^n_{\widetilde{x}}$ determines $\widetilde{x}$
uniquely because $\B=C^*(
\bigcup_{n=0}^\infty\delta_*^n ({\cal A}))$. Since $\delta_*$ is
an endomorphism of $\B$
 the functionals $\xi^n_{\widetilde{x}}$ are linear and multiplicative
on $\A$. So either $\xi^n_{\widetilde{x}}\in X$ ($X$ is spectrum
of $\A$) or $\xi^n_{\widetilde{x}} \equiv 0$. It follows then that
the mapping
\begin{equation}
\label{xi000'} M(\B) \ni {\widetilde{x}} \to
(\xi^0_{\widetilde{x}}, \xi^1_{\widetilde{x}}, ...)\in
\prod_{n=0}^{\infty} (X\cup\{0\})
\end{equation} is an injection and  the following statement is true, see Theorems 3.1 and  3.3 in
\cite{maxid}.
\begin{thm}\label{ideals}
 Let
$\delta(\cdot)=U(\cdot )U^*$ be an endomorphism of ${\cal A}$, $
  U^*U\in \A'$, and $\alpha:\Delta_1\rightarrow X$ be the partial mapping
  determined by $\delta$.
Then the mapping (\ref{xi000'}) defines a topological embedding
of $M(\B)$ into topological
space $\bigcup_{N=0}^{\infty}\widehat{X}_N\cup X_\infty$. Under this
embedding we have
$$\label{equality}
\bigcup_{N=0}^{\infty}X_N\cup X_\infty \subset
M(\B)\subset \bigcup_{N=0}^{\infty}\widehat{X}_N\cup
X_\infty
$$
where
$$
\widehat{X}_N=\{\widetilde{x}=(x_0,x_1,...,x_N,0,...): x_n\in \Delta_n,\,
\alpha(x_{n})=x_{n-1},\,1\leq n\leq N\},$$
$$
X_\infty=\{\widetilde{x}=(x_0,x_1,...): x_n\in \Delta_n,\,
\alpha(x_{n})=x_{n-1},\, \,1\leq n\}.
$$
$$
X_N=\{\widetilde{x}=(x_0,x_1,...,x_N,0,...)\in
\widehat{X}_N:x_N\notin \Delta_{-1}\},
$$
The  topology on
$\bigcup_{N=0}^{\infty}\widehat{X}_N\cup X_\infty$ is defined  by a
fundamental system of neighborhoods of points $\widetilde{x}\in
\widehat{X}_N$ given by
$$
 O(a_1,...,a_k,\varepsilon)=\{\widetilde{y}\in \widehat{X}_N:
|a_i(x_N)-a_i(y_N)|<\varepsilon, \, \, i=1,...,k\}
$$
and respectively of   $\widetilde{x}\in X_\infty$ by
$$ O(a_1,...,a_k,n, \varepsilon)=\{\widetilde{y}\in
\bigcup_{N=n}^{\infty}\widehat{X}_N\cup X_\infty:
|a_i(x_n)-a_i(y_n)|<\varepsilon, \, \, i=1,...,k\}
$$
 where $\varepsilon > 0$,
$a_i\in \A$ and $k,n \in  \N$.
\end{thm}
\begin{rem}\label{remark2.3}
The topology on $X$ is weak$^*$. One immediately sees then (see
(\ref{xi000'})), that the topology on $\bigcup_{N\in \N}\widehat{X}_N\cup
X_\infty$ is in fact the product topology inherited from
$\prod_{n=0}^{\infty} (X\cup\{0\})$ where $\{0\}$ is clopen.
\end{rem}
\noindent The foregoing theorem gives us an estimate of
$M(\B)$ and aiming at sharpening  that
result we need to strengthen the assumptions. If we replace the condition $U^*U\in\A'$ with the stronger one
\begin{equation}\label{condition}
 U^*U\in \A,
 \end{equation}
then the full information on $\B$ is carried
by the pair $(\A,\delta)$, cf. \cite[Theorem 3.4]{maxid}.
 \begin{thm}\label{ideals2}
Under the assumptions of Theorem \ref{ideals} with $U^*U\in\A'$
replaced by $U^*U\in \A$ we get
$$
M(\B)=\bigcup_{N=0}^{\infty} X_N\cup X_\infty.
$$
 \end{thm}

 This motivates us to take a closer  look at  condition  (\ref{condition}).\\
 Firstly, let us observe
 \cite[3.5]{maxid} that if $U^*U\in\A'$ then $\delta$ is an endomorphism of the
$C^*$-algebra ${\cal A}_1=C^*({\cal A},U^*U )$  and since
$C^*(\bigcup_{n=0}^\infty U^{*n}\A_1 U^n)=\B$ one can
apply the preceding theorem to the algebra $\A_1$  for the full
description of $M(\B)$.
This procedure turns out to be very fruitful in many situations.

\begin{exm}\label{maxidexm}
Let the elements of $\A$ be the   operators   of multiplication by
periodic sequences of period $n$, on the Hilbert space $l^2(\N)$,
and let   $U$ be the  co-isometry given by $ [Ux](k)= x(k+1),$ for  $x\in
l^2(\mathbb{N})$, $k\in \N$.
 Then $X=\{x_0,...,x_{n-1}\}$ and $\al(x_k)=x_{k+1\,\,(mod\,\, n)}$. If for $k=0,...,n-1$ we write
 $$
(\infty,k)=(x_k,x_{k-1},...,x_1,x_0,x_{n-1},x_{n-2},\dots)$$
and
$$(N,k)=(\underbrace{x_k,x_{k-1},...,x_1,x_0,x_{n-1},...}_N,x_{n-r},0,0,\dots)
$$
where $N-r\equiv k$ $(mod\,\, n)$ (for each $N$ there are $n$ pairs $(N,k)$),
 then from Theorem
\ref{ideals} we get
$$
 \{\infty\}\times \{0,1,...,n-1\}\subset M(\B) \subset \overline{\N}\times \{0,1,...,n-1\}
$$
where $\overline{\N}=\N\cup\{\infty\}$ is a compactification
 of the discrete space $\N$. In order  to describe $M(\B)$  precisely let us pass to
 the algebra  ${\cal A}_1=C^*({\cal A},U^*U )$.
As $U^*U$ is the operator of multiplication by $(0,1,1,...)$ the
elements of $\A_1$ are the operators   of multiplication by sequences of the form
$(a,h(0),h(1),...)$ where $a$ is arbitrary and $h(k+n)=h(k)$, for
all $k\in \N$. In an obvious manner (with a slight abuse of
notation) we infer the spectrum of $\A_1$ to be
$\{y,x_0,...,x_{n-1}\}$, and the mapping generated by $\delta$
considered as an endomorphism of $\A_1$  acts as follows:
$\alpha(x_k)=x_{k+1\,\,(mod\,\, n)}$ and $\alpha(y)=x_0$. Abusing notation once again and putting $$(N,k)=(\underbrace{x_k,x_{k-1},...,x_1,x_0,x_{n-1},...,x_0}_N,y,0,0,...)$$ where $N\equiv k$ $(mod\,\, n)$ (for each $N$ there is now the only one pair $(N,k)$), in view
of Theorem  \ref{ideals2}  we have
$$
M(\B)=\{ (N,k)\in \N\times \{0,1,...,n-1\}:
N\equiv k \,\,(mod\,\,n)\} \cup \big(\{\infty\}\times \{0,1,...,n-1\}\big),
$$
so $M(\B)$ can be imagined as a  spiral subset
of the  cylinder $\overline{\N}\times \{0,1,...,n-1\}$, (see
Figure 1).
\end{exm}

\setlength{\unitlength}{.8mm}
\begin{picture}(120,78)(-30,-4)

%podstawa cylindra

 \put(35,10){\circle*{1}}
\put(28,6){\scriptsize{$(\infty,0)$}} \put(55,16){\circle*{1}}
\put(51,12){\scriptsize{$(\infty,n\,\,\,1)$}}
\put(59.7,12.5){\tiny{$-$}}

\put(53,4){\circle*{1}} \put(49,0.5){\scriptsize{$(\infty,1)$}}
\put(75,4){\circle*{1}} \put(71,0.5){\scriptsize{$(\infty,2)$}}

\put(53,16){\line(-3,-1){16}} \put(37,9){\line(3,-1){14}}
\qbezier[60](76,16)(66,16)(56,16) \put(55,4){\line(1,0){18}}

%kropki

%\put(78,5){.} \put(81,6){.} \put(84,7.25){.}
\put(88.5,9.2){.} \qbezier[40](76,4.5)(82,6.5)(88,9)

\put(85.7,10.7){.} \put(83,12.5){.} \put(80,14.2){.}
\put(76.55,15.7){.}

% pionowe kreski

\qbezier[100](35,12)(35,40)(35,68)
\qbezier[40](55,18)(55,43)(55,74)
\qbezier[100](53,6)(53,35)(53,62)
\qbezier[100](75,6)(75,35)(75,62)
\qbezier[40](89.5,11)(89.5,40)(89.5,68)
\qbezier[40](77,18)(77,43)(77,73)

% punkty

\put(35,69){\circle*{1}}\put(26.7,66.2){\scriptsize{$(0,0)$}}
\put(53,59){\circle*{1}} \put(44.7,56){\scriptsize{$(1,1)$}}
\put(75,55.7){\circle*{1}} \put(66.7,53){\scriptsize{$(2,2)$}}
 \put(89,55.5){.}
 %\put(76.55,55.5){.}
 \put(54.55,54.5){.}

\put(35,48){\circle*{1}}\put(26.7,45){\scriptsize{$(n,0)$}}
\put(53,38){\circle*{1}}
\put(41,35){\scriptsize{$(n\,\,\,1,1)$}}\put(44.5,35.5){\tiny{$+$}}

\put(75,34.7){\circle*{1}}
\put(63,32){\scriptsize{$(n\,\,\,2,2)$}}\put(66.5,32.5){\tiny{$+$}}

\put(89,34.5){.}

 % linie
 \qbezier(36,68.5)(52,59.5)(52,59.5)
\qbezier(54,59)(54,59)(74,55.7)
\qbezier[32](76,55.7)(81.5,55.7)(88.5,55.7)
\qbezier[30](54,54.5)(45,51.25)(36,48)
\qbezier(36,47.5)(52,38.5)(52,38.5)
\qbezier(54,38)(54,38)(73.7,34.8)
\qbezier[32](76,34.7)(81.5,34.7)(88.5,34.7)

%Podpis
\put(110, 42){\textbf{Figure 1.}}

 \put(110, 34){Maximal ideal space of}

\put(110, 28){the coefficient algebra}

 \put(110, 22){from Example \ref{maxidexm}.}
\end{picture}
\par
  Condition  (\ref{condition}) is closely related
to the openness of $\Delta_{-1}$ (as $\Delta_1$ is compact and
$\al$ is continuous $\Delta_{-1}$ is always closed).
\begin{prop}\label{prop(O)sition}
Let $P_{\Delta_{-1}}\in \A$ be the projection corresponding to the
characteristic function $\chi_{\Delta_{-1}}\in C(X)$. If
$U^*U\in\A$ then $\Delta_{-1}$ is open and $U^*U =
P_{\Delta_{-1}}$. If $U^*U\in\A'$, $\Delta_{-1}$ is open and $\A$
acts nondegenerately on $H$, then $
 U^*U \leqslant P_{\Delta_{-1}}. $
\end{prop}
\begin{Proof}
Let $U^*U\in\A$.  We show that the image $\Delta_{-1}$ of the set
$\Delta_1$ of functionals satisfying (\ref{e1.0}) under the
mapping (\ref{e2.0}) is the set of functionals $x\in
X$ satisfying $x(U^*U)=1$.\\
Let $x'\in \Delta_1$, that is $x'(UU^*)=1$. Putting $x=\al(x')$ we
have $ x(U^*U)=x'(\delta(U^*U))=x'(UU^*UU^*)=x'(UU^*)x'(UU^*)=1. $
Now, let $x\in X$ be such that $x(U^*U)=1$. We define on
$\delta(\A)$ a multiplicative functional $x'(b):=x(U^*bU)$,
$b\in\delta(\A)$. For $b=\delta(a)$, $a\in \A$, we then have
$$
x'(\delta(a))=x(U^*UaU^*U)=x(U^*U)x(a)x(U^*U)=x(a).
$$
 $x'$ is therefore well defined. Since
$x'(\delta(1))=x(1)=1$ it is  non zero and there exists  its
extension $\overline{x}'\in X$ on $\A$ (see \cite[2.10.2]{Dix}).
It is clear that $\overline{x}'\in \Delta_{1}$ and
$\alpha(\overline{x}')=x$. Hence $x\in \Delta_{-1}$. Thus we have
proved that $ x\in \Delta_{-1}\Longleftrightarrow x(U^*U)=1 $,
which means that $U^*U\in \A$ is the characteristic
function of $\Delta_{-1}$.
It follows then that  $\Delta_{-1}$ is clopen.\\
Now, let $\Delta_{-1}$ be open. Then $\chi_{\Delta_{-1}}\in C(X)$
and
$\delta(\chi_{\Delta_{-1}})=\chi_{\al^{-1}(\Delta_{-1})}=\chi_{\Delta_{1}}$.
As $\A$ acts nondegenerately, rewriting this equation in terms of
operators we have
\begin{equation}\label{UPU}
 UP_{\Delta_{-1}}U^*=UU^*.
 \end{equation}
 Letting $H_i=U^*UH$
be the initial and $H_f=UU^*H$  be the final space of $U$ we get
$U^*:H_f\rightarrow H_i$ is an isomorphism and $U:H_i\rightarrow
H_f$ is its inverse. Taking arbitrary $h\in H_f$ and applying the both sides of (\ref{UPU}) to it
 we obtain $ UP_{\Delta_{-1}}U^*h=h
$, and hence $P_{\Delta_{-1}}H_i=H_i$, that is $U^*U \leqslant
P_{\Delta_{-1}}$.
\end{Proof}
\textbf{Note.} The inequality in the second part of the preceding
proposition can not be replaced by
  equality. In order to see that consider for instance $\A$ and $U$ from Example \ref{maxidexm}.
\\
By virtue of Proposition \ref{prop1} the mappings $\delta$ and
$\delta_*$ are  endomorphisms of the $C^*$-algebra
$\B$.
 With the help of the presented
theorems we can now find the form of the partial mappings they
generate. We shall rely on the fact \cite[2.5]{maxid} expressed
by the coming proposition.
\begin{prop}
\label{tu1} Let  $\delta(\cdot)=U(\cdot) U^*$ and
$\delta_*(\cdot)=U^*(\cdot)U$ be  endomorphisms of $\A$ and let
$\alpha$ be the partial mapping of $X$ generated by $\delta$.
 Then   $\Delta_1$
and   $\Delta_{-1}$  are clopen  and $\alpha:\Delta_1 \rightarrow
\Delta_{-1}$ is a homeomorphism. Moreover, the endomorphism $\delta_*$
is given on $C(X)$ by the formula \begin{equation}
 \delta_*(a)(x)=\left \{ \begin{array}{ll} a(\alpha^{-1}(x))& ,\
\ x\in
\Delta_{-1}\\
0 & ,\ \ x\notin \Delta_{-1} \end{array} \right. \end{equation}
\end{prop}
Finally we arrive at the closing theorem.
\begin{thm}\label{tu3}
Let the hypotheses  of Theorem \ref{ideals} hold. Then
\begin{description}
\item[i)] the sets
$$
\widetilde{\Delta}_1=\{(x_0,...)\in
M(\B): x_0 \in \Delta_1\},$$
$$
\widetilde{\Delta}_{-1}=\{(x_0,x_1,...)\in
M(\B): x_1\neq 0 \}
$$ are clopen subsets of
$M(\B)$, \item[ii)] the endomorphism $\delta$ generates
on $M(\B)$ the partial homeomorphism
$\widetilde{\alpha}:\widetilde{\Delta}_1 \rightarrow
\widetilde{\Delta}_{-1}$  given by the formula
\begin{equation}\label{alfazfalkom}
\widetilde{\alpha}(x_0,...)=(\alpha(x_0),x_0,...),\qquad
(x_0,...)\in \widetilde{\Delta}_1,
\end{equation}
 \item[iii)]
the partial mapping generated by  $\delta_*$ is the inverse of
$\widetilde{\alpha}$, that is
$\widetilde{\alpha}^{-1}:\widetilde{\Delta}_{-1}\rightarrow
\widetilde{\Delta}_{1}$ where
\begin{equation}\label{alfa-1zfal}
\widetilde{\alpha}^{-1}(x_0,x_1...)=(x_1,...),\qquad
(x_0,x_1...)\in \widetilde{\Delta}_{-1}. \end{equation}
\end{description}
\end{thm}
\begin{Proof}  We rewrite Proposition \ref{tu1} in terms of Theorem  \ref{ideals}.
Let $\widetilde{x}=(x_0,x_1,...)\in M(\B)$. From
(\ref{e1.0}) we get
$$\widetilde{x}\in \widetilde{\Delta}_1 \Longleftrightarrow \widetilde{x}(UU^*)=1,\qquad
\widetilde{x}\in \widetilde{\Delta}_{-1} \Longleftrightarrow
\widetilde{x}(U^*U)=1.$$ However, the definition (\ref{xi2}) of
functionals $\xi_{\widetilde{x}}^n=x_n$ implies that
$\widetilde{x}(UU^*)=1\Longleftrightarrow x_0(UU^*)=1,$  and
$\widetilde{x}(U^*U)=1\Longleftrightarrow x_1(1)=1,$
which proves  i).\\
The mapping $\widetilde{\alpha}$  generated by
$\delta$ on  $M(\B)$  (see  (\ref{e2.0})), is given
by the composition: $ \widetilde{\alpha}(\widetilde{x})\equiv
\widetilde{x}\circ \delta. $ So, let
$\widetilde{x}=(x_0,x_1,...)\in \widetilde{\Delta}_1$, then the
sequence of functionals  $\xi^n_{\widetilde{x}}\,$ satisfies:
$\xi^n_{\widetilde{x}}(a)=a(x_n)$, $a\in \A$, $n\in\N$. Now let us
consider an analogous sequence of functionals
$\xi^n_{\widetilde{\alpha}(\widetilde{x})}\,$ defining the point
$\widetilde{\alpha}(\widetilde{x})=(\overline{x}_0,\overline{x}_1,...)$.
For $n>0$ we have
$$a(\overline{x}_n)=\xi^n_{\widetilde{\alpha}(\widetilde{x})}(a)=\widetilde{\alpha}(\widetilde{x})(\delta_*^n(a))=
\widetilde{x}(\delta(\delta_*^n(a)))=\widetilde{x}(UU^{*n}aU^nU^*)=$$
$$=\widetilde{x}(UU^*)\widetilde{x}(\delta_*^{n-1}(a))\widetilde{x}(U^*U)=
\widetilde{x}(\delta_*^{n-1}(a))=\xi^{n-1}_{\widetilde{x}}(a)=a(x_{n-1}),$$
while  for $n=0$ we have
$$
a(\overline{x}_0)=\xi^0_{\widetilde{\alpha}(\widetilde{x})}(a)=\widetilde{\alpha}(\widetilde{x})(a)
=\widetilde{x}(\delta(a))=
\xi^0_{\widetilde{x}}(\delta(a))=\delta(a)(x_0)=a(\alpha(x_0)).$$
Thus we infer that
$\widetilde{\alpha}(\widetilde{x})=(\alpha(x_0),x_0,...).$ By
Proposition  \ref{tu1}, $\widetilde{\alpha}^{-1}$ is the inverse
to mapping $\widetilde{\alpha}$. Hence we  get (\ref{alfa-1zfal}).
\end{Proof}

\section{Reversible  extension of a partial dynamical system}
One of the most important consequences of Theorems
 \ref{ideals2} and \ref{tu3} is that although the algebra  $\B$
 is relatively bigger than $\A$ and its structure depends  on $U$ and $U^*$
 ($U$ and $U^*$ need not  be in $\A$) the $C^*$-dynamical system $(\B,\delta)$
 still can be reconstructed
 by means of the intrinsic features of $(\A,\delta)$ itself (provided (\ref{condition}) holds).
  Therefore in this section,
  we make an effort to investigate effectively this reconstruction and, as it is purely topological,
   we forget for the time being about its algebraic aspects.
   \\
    Once having the system $(X,
    \al)$, we will construct a pair $(\X,\tal)$:
    \begin{itemize}
   \item a compact space $\X$ - a counterpart
   of the maximal ideal space  obtained in Theorem \ref{ideals2}
    (or  a 'lower estimate' of it, see Theorem \ref{ideals}),
    and
    \item a partial injective mapping $\tal$ - a counterpart of the mapping from Theorem
    \ref{tu3}.
    \end{itemize}
     We then show  some useful results about the structure
    of $(\X,\tal)$. In particular, we calculate $(\X,\tal)$ for
    topological Markov chains, and  show the interrelation between $
    \X$ and  projective limits.
\\
The most interesting case occurs  when the
 injective  mapping $\tal$ has an open image. We shall call   such mappings
 partial homeomorphisms, cf \cite{exel_laca_quigg}. More
precisely,
\\

a \emph{partial homeomorphism} is  a partial  mapping  which is
injective and has open image.
\\

\noindent Let us recall that by a partial mapping of  $X$ we always
mean a continuous mapping $\al:\Delta_1\rightarrow X$ such that
$\Delta_1\subset X$ is clopen, and so if $\al$ is  a partial
homeomorphism, then $\al^{-1}:\Delta_{-1}\rightarrow X$ is a
partial mapping of $X$ in our sense, that is $\Delta_{-1}$ is
clopen.
\\
 As we shall see if $\Delta_{-1}$ is open then $\tal$
is a partial homeomorphism and $\X$ actually becomes the spectrum of
a certain coefficient algebra (see Theorem \ref{coefalgebra}). That is the reason why we shall
often assume the openness of $\Delta_{-1}$. However in this section we do not make it a standing
assumption in order to get to know better the role
of it and the condition (\ref{condition}), cf. Proposition
\ref{prop(O)sition}.
\subsection{The system $(\X,\tal)$}
Let us fix a partial dynamical system $(X, \al)$ and let us
consider a disjoint union $X\cup\{0\}$ of the set $X$ and the
singleton $\{0\}$ (we treat here $0$  as a symbol rather than the
number). We define   $\{0\}$ to be  clopen
 and hence $X\cup\{0\}$ is a compact topological space. We will let $\X$ to be  a subset
$$
\X\subset \prod_{n=0}^{\infty} (X\cup\{0\})
$$
of the product of $\aleph_0$ copies of the space $X\cup\{0\}$
where the elements of $\X$ represent anti-orbits of the partial
mapping $\al$. Namely we set \begin{equation}\label{widetildeX}
\X=\bigcup_{N=0}^{\infty}X_N\cup X_\infty \end{equation} where
$$
X_N=\{\x=(x_0,x_1,...,x_N,0,...): x_n\in \Delta_n,\,x_N\notin
\Delta_{-1},\, \al(x_{n})=x_{n-1},\,\, n=1,...,N\},
$$
$$
X_\infty=\{\x=(x_0,x_1,...): x_n\in \Delta_n,\,
\al(x_{n})=x_{n-1},\, \,n\geqslant 1\}.
$$
 The natural topology on  $\X$ is the one induced from the space
$\prod_{n=0}^{\infty} (X\cup\{0\})$ equipped with the product
topology, cf. Remark \ref{remark2.3}. Since $ \X\subset
\prod_{n=0}^{\infty} (\Delta_n\cup\{0\}), $
 the
topology on $\X$ can also be regarded as the topology inherited
from $\prod_{n=0}^{\infty} (\Delta_n\cup\{0\})$.
\begin{defn}\label{widetildeXdefn}
We shall call the topological space $\X$  the \emph{extension of
$X$ under $\al$}, or briefly the \emph{ $\al$-extension} of $X$.
\end{defn}

\begin{thm}\label{compactnes} The subset $X_\infty$ is  compact  and the subsets $X_N$ are clopen in $\X$.
 Moreover, the following conditions are equivalent:
\\
a) $\Delta_{-1}$ is open.\\
b)  $\X$ is   compact.\\
c)  $X_0$ is compact.\\
d) $X_N$ is compact for every $N\in \N$.
\end{thm}
\begin{Proof} As the  sets $\Delta_n$, $n\in \N$,
are clopen, by Tichonov's theorem, the space $\prod_{n=0}^{\infty}
(\Delta_n\cup\{0\})$ is compact,  and
 to prove the compactness  of $X_\infty$
it suffices to show that $X_\infty$ is a closed subset of
$\prod_{n=0}^{\infty} (\Delta_n\cup\{0\})$, or equivalently that its complement is open. To this end, let
$\x=(x_0,x_1,...)\in \prod_{n=0}^{\infty} (\Delta_n\cup\{0\})$ and suppose that $\x\notin
X_\infty$.  We show that $\x$ has an open neighborhood contained in
 $\prod_{n=0}^{\infty} (\Delta_n\cup\{0\})\setminus \X_\infty$.
\\
In view of the definition of $X_\infty$ we have two possibilities:
\Item{1)}  there is $n>0$ such that  $x_n=0$,
\Item{2)} there is $n>0$ such that  $x_n\in \Delta_n$,  and
$\al(x_n)\neq x_{n-1}$,
\\
If  1) holds then  the set
$\V=\prod_{k=0}^{\infty}V_k$ where $V_k=\Delta_k\cup\{0\}$, for
$k\neq n$, and $V_n=\{0\}$,  is an open
neighborhood of $\x$ and $\V\cap X_\infty=\emptyset$. \\
Let us suppose now that 2) holds. We may also suppose that  $x_{n-1}\neq 0$ $(x_{n-1}\in\Delta_{n-1})$. Hence there
exist two disjoint open subsets $V_1$, $V_2\subset \Delta_{n-1}$
such that $\al(x_n)\in V_1$ and $x_{n-1}\in V_2$. Clearly, the set
$\V=\prod_{k=0}^{\infty}V_k$ where $V_k=\Delta_k\cup\{0\}$, for
$k\neq n-1,n$, and $V_{n-1}=V_2$,  $V_n=\al^{-1}(V_1)$,
 is an open
neighborhood of $\x$, and  $V_1\cap V_2=\emptyset$ guarantees that  $\V\cap X_\infty=\emptyset$.
\\
Fix $N\in \mathbb{N}$. To prove that $X_N$  is  open we  recall that
$\Delta_n$, $n\in \N$, are clopen and $\Delta_{-1}$ is closed.
Hence $\Delta_n\setminus\Delta_{-1}$, $n\in \N$,  are open, and it is easy to see
that
$$
X_N=\X\cap (\Delta_0\times \Delta_1\times ...\times
\Delta_{N-1}\times \Delta_N\setminus\Delta_{-1}
\times(\Delta_{N+1}\cup\{0\})\times ...).
$$
Hence  $X_N$ is an open subset of $\X$. It also closed because its complement is a sum of two open sets:
$\V_1=\{(x_0,x_1...)\in \X: x_N=0\}$ and $\V_2=\{(x_0,x_1...)\in \X: x_{N+1}\neq 0\}$.
\\
We prove now the equivalence of a), b), c) and d).
\\
a)$\Rightarrow$b). Suppose that $\Delta_{-1}$ is open. We
  prove the  compactness  of $\X$ in an analogous fashion as we proved the compactness of $X_\infty$.
Let
$\x=(x_0,x_1,...)\in \prod_{n=0}^{\infty} (\Delta_n\cup\{0\})\setminus \X$. In view of the definition of $\X$
 we have the  three possibilities:

\Item{1)} there are $n,m\in \N$  such that  $x_n=0$ and $x_{n+m}\in
\Delta_{n+m}$ ($x_{n+m}\neq 0$),
\Item{2)} there is $n>0$ such that  $x_n\in \Delta_n$,  and
$\al(x_n)\neq x_{n-1}$,
\Item{3)} for some $n>0$ we have  $x_n\in \Delta_n\cap \Delta_{-1}$,  and
$x_{n+1}=0$.
\\
Let us suppose that 1) holds.  Then the set
$\V=\prod_{k=0}^{\infty}V_k$ where $V_k=\Delta_k\cup\{0\}$, for
$k\neq n,n+m$, and $V_n=\{0\}$, $V_{n+m}=\Delta_{n+m}$, is an open
neighborhood of $\x$.  It is clear that none of the points from
$\X$
belongs to $\V$. \\
The same  argumentation  as the one concerning $X_\infty$ shows that in the case
 2) $\x$ lies in the interior of $\prod_{n=0}^{\infty} (\Delta_n\cup\{0\})\setminus \X$.
\\
If we suppose that 3) holds, then the set
$\V=\prod_{k=0}^{\infty}V_k$ where $V_k=\Delta_k\cup\{0\}$, for
$k\neq n,n+1$, and $V_n=\Delta_n\cap \Delta_{-1}$, $V_{n+1}=\{0\}$, is an open
neighborhood of $\x$ (here we use the openness of $\Delta_{-1}$). Clearly $\V$ does not contain any point from
$\X$.
\\
b)$\Rightarrow$c).  $X_0$ is closed and hence it is compact.\\
c)$\Rightarrow$a).
Suppose that $\Delta_{-1}$ is not open. Then $X\setminus \Delta_{-1}$ is not closed and hence it is  not compact.
Thus there is an  open cover $\{V_i\}_{i\in I}$  of $X\setminus \Delta_{-1}$ which does not  admit
 a finite subcover. Defining $\V_i=\{(x_0,x_1...)\in \X: x_0\in V_i\}$, for $i\in I$, we get an open cover
  of $X_0$ which does not  admit  a finite subcover. Hence $X_0$ is not compact.
\\ Thus we see that a)$\Leftrightarrow$b)$\Leftrightarrow$c). As b)$\Rightarrow$d) and d)$\Rightarrow$c) are obvious,
the proof is complete.
 \end{Proof}
It is interesting how $\X$ depends on $\al$. For instance,  if
$\al$ is surjective then $X_N$, $n\in \N$, are empty and
$\X=X_\infty$, in this case $\X$ can be defined as a projective
limit, see Proposition \ref{projectivelimit}. If  $\al$ is
injective then  a natural continuous projection $\Phi$ of $\X$
onto $X$ given by the  formula
\begin{equation} \label{Phimap} \Phi(x_0,x_1,...)=x_0
\end{equation}
 becomes a bijection and, as we will see, in case $\Delta_{-1}$ is open even a homeomorphism.
  But the farther from injectivity $\al$ is,
the farther $\X$ is from $X$.
\begin{prop}\label{Phihomeomorphism}
Let $(X,\al)$ be such a  system that $\al$ is injective. Then $\Phi:\X\rightarrow X$ is a
homeomorphism if and only if
$\Delta_{-1}$ is open.
\end{prop}
\begin{Proof} If $\Delta_{-1}$ is not open then by Theorem \ref{compactnes},
$\X$ is not compact and hence not homeomorphic to $X$.
Suppose then that $\Delta_{-1}$ is  open. Then $\al:\Delta_1\rightarrow X$  is an open mapping because $\al$ is a homeomorphism of compact set $\Delta_1$ onto the compact set $\Delta_{-1}$.
 We only need to show that the mapping
$\Phi^{-1}$ is continuous or, which is the same, that $\Phi$ is
open. To see this it is enough to  look at a subbase for the
topology in $\X$ in a appropriate way. Indeed, let $\U\subset \X$
be of the   form
$$\U=\X\cap (U_0\times U_1\times ...\times
U_{N}\times (\Delta_{N+1}\cup\{0\})
\times(\Delta_{N+2}\cup\{0\})\times ...)
$$
where $U_i\subset X\cup\{0\}$, $i=1,...,N$, are open. Without loss
of generality we can assume that $U_i\subset \Delta_i\cup\{0\}$,
$i=1,...,N$. There are the two possibilities:
 \Item{1)} If $0\notin U_N$ then the set $U_{N-1}\cap\al(U_N)$ is
open and
$$
\U=\X\cap (U_0\times U_1\times ...\times
\Big(U_{N-1}\cap\al(U_N)\Big) \times (\Delta_{N}\cup\{0\})
\times(\Delta_{N+1}\cup\{0\})\times ...)
$$
\Item{2)} If $0\in U_N$ then the set
$U_{N-1}\cap\al(U_N\setminus\{0\})\cup U_{N-1}\setminus
\Delta_{-1}\cup \{0\}$ is open and
$$
\U=\X\cap (U_0\times U_1\times ...\times
\Big(U_{N-1}\cap\al(U_N\setminus\{0\})\cup
U_{N-1}\setminus\Delta_{-1}\cup\{0\}\Big) \times
(\Delta_{N}\cup\{0\}) \times ...).
$$
Applying the above procedure   $N$ times  we conclude that
$\U=\X\cap \prod_{k=0}^{\infty}V_k$ where $V_k=\Delta_k\cup\{0\}$,
for $k>0$, and $V_0\subset X\cup\{0\}$ is a certain  open set. It
is obvious that $\Phi(\U)=V_0$. Thus $\Phi$ is an open mapping and
the proof is complete.
\end{Proof}
\begin{exm}\label{q-exm} Consider the dynamical system $(X,\al)$ where $X=[0,1]$ and
$\al(x)=q\cdot x$, for fixed $0<q<1$ and any $x\in[0,1]$. In this
situation  $\al$ is injective and $\Delta_{-1}=[0,q]$ is not open.
As $\Phi$ is bijective, we can identify $\X$ with an interval
$[0,1]$ but the topology on $\X$ differs from the natural topology
on $[0,1]$. It is not hard to check that the topology on
$\X=[0,1]$ is generated by intervals $[0,a)$, $(a,b)$, $(b,1]$,
where $0<a<b<1$, and singletons $\{q^k\}$, $k > 0$. Thus $\X$ is not compact and $\Phi$ is
not a homeomorphism.
\end{exm}

Now, we would like to investigate a partial mapping $\tal$ on $\X$
associated with $\al$. It seems very natural to look for a partial
mapping $\tal$ such that $ \Phi\circ \tal=\al \circ \Phi$
 wherever the superposition $\al \circ \Phi$ makes sense.
  If such  $\tal$ exists then its  domain is  $\TDelta_1:=\Phi^{-1}(\Delta_{1})$
  and  its image is contained in
 $\TDelta_{-1}:=\Phi^{-1}(\Delta_{-1})$. Moreover, as $\Phi$ is continuous, we get
$$
\TDelta_{1}=\Phi^{-1}(\Delta_1)=\{\x=(x_0,x_1,...)\in \X: x_0 \in
\Delta_{1}\}
$$
is clopen while
$$
\TDelta_{-1}=\Phi^{-1}(\Delta_{-1})=\{\x=(x_0,x_1,...)\in \X: x_0
\in \Delta_{-1}\}
$$
is closed  and if $\Delta_{-1}$ is open then $\TDelta_{-1}$ is
open too. It is not a surprise that there  always exists  a
homeomorphism $\tal$ such that the following diagram
\begin{equation}\label{phidiagram}
 \begin{array}{c c c} \TDelta_1 & \stackrel{\tal}{\longrightarrow} &\TDelta_{-1}  \\
\Phi\downarrow &      & \downarrow \Phi\\
\Delta_1 & \stackrel{\al}{\longrightarrow} & \Delta_{-1}
\end{array}\end{equation}
commutes. However, the commutativity of the diagram
(\ref{phidiagram}) does not determine the homeomorphism $\tal$
uniquely.
\begin{prop}\label{proptildealpha}
The mapping $\tal:\TDelta_1\rightarrow \TDelta_{-1}$ given by the
formula \begin{equation}\label{tildealpha}
\tal(\x)=(\al(x_0),x_0,...),\qquad \x=(x_0,...)\in \TDelta_1
\end{equation}
 is  a  homeomorphism (and hence if $\Delta_{-1}$ is open $\tal$ is a partial homeomorphism of $\X$) such that  the
diagram (\ref{phidiagram}) is commutative.
 \end{prop}
\begin{Proof} The inverse of $\al$ is given by the formula $\tal^{-1}((x_0,x_1,...))=(x_1,x_2,...)$.
 The straightforward  equations $\tal(\TDelta_1\cap (U_0\times U_1 \times
...))=\X\cap (X\times U_0\times ...)$  and
$\tal^{-1}(\TDelta_{-1}\cap (U_0\times U_1\times
U_2\times...))=\X\cap ((\al^{-1}(U_0)\cap U_1)\times U_2
\times...)$, and  the definition of the topology in $\X$  imply that $\tal$ and $\tal^{-1}$ are continuous and
hence they are homeomorphisms.  The commutativity of the diagram
(\ref{phidiagram}) is obvious.
\end{Proof}

In this manner we can attach to every pair $(X,\al)$ such that
$\Delta_{-1}$ is open another system $(\X,\tal)$ where $\tal$ is a
partial homeomorphism of $\X$ and if $\al$ is a partial
homeomorphism   then systems $(X,\al)$ and $(\X,\tal)$ are
topologically equivalent via   $\Phi$, see Proposition
\ref{Phihomeomorphism}. In particular case  of a classical
irreversible dynamical system, that is when $\al$ is a covering
mapping of $X$,  $\tal$ is a full homeomorphism and hence
$(\X,\tal)$ is a classical reversible dynamical system. This
motivates us to coin the following definition.
\begin{defn} Let $(X,\al)$ be  a partial dynamical system and let $\Delta_{-1}$ be open.  We say that the pair
$(\X,\tal)$, where $\X$ and $\tal$ are given by (\ref{widetildeX})
and (\ref{tildealpha}) respectively, is a \emph{reversible
extension} of  $(X,\al)$.
\end{defn}
\begin{exm}\label{simplexample} It may happen that two different partial dynamical systems
have the same reversible extension. Let $(X,\al)$  and $(X',\al')$ be given by the relations:
 $X=\{x_0,x_1,x_2,y_2\}$, $\Delta_1=X\setminus\{x_0\}$,  $\al(y_2)=\al(x_2)=x_1$,
   $\al(x_1)=x_0$; $X'=\{x_0',x_1',x_2',y_1',y_2'\}$, $\Delta_1'=X'\setminus\{x_0'\}$
    and $\al'(x_2')=x_1'$, $\al'(y_2')=y_1'$ ,  $\al'(y_1')=\al'(x_1')=x_0'$; or by the diagrams:\\
\setlength{\unitlength}{.8mm}
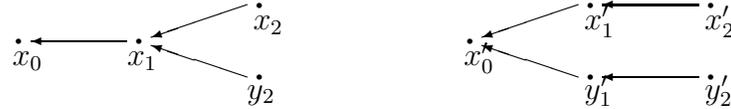
\begin{picture}(150,20)(-40,0)
\put(15,10){\circle*{1}} \put(14,6){$x_0$}
\put(35,10){\circle*{1}} \put(33,6){$x_1$}
\put(55,16){\circle*{1}} \put(54,12){$x_2$}
\put(55,4){\circle*{1}} \put(53,0.5){$y_2$}
\put(53,16){\vector(-3,-1){16}} \put(33,10){\vector(-1,0){16}}
\put(53,4){\vector(-3,1){16}}

\put(90,10){\circle*{1}} \put(89,6){$x_0'$}
\put(110,16){\circle*{1}} \put(109,12){$x_1'$}
\put(110,4){\circle*{1}} \put(109,0.5){$y_1'$}
\put(130,16){\circle*{1}} \put(129,12){$x_2'$}
\put(130,4){\circle*{1}} \put(129,0.5){$y_2'$}

\put(108,16){\vector(-3,-1){16}} \put(108,4){\vector(-3,1){16}}
\put(128,16){\vector(-1,0){16}} \put(128,4){\vector(-1,0){16}}

\end{picture}\\
Then   $\X$ consist of the following points
$\x_0=(x_0,x_1,x_2,0,...)$, $\x_1=(x_1,x_2,0,...)$,
$\x_2=(x_2,0,...)$,  $\y_0=(x_0,x_1,y_2,0,...)$,
$\y_1=(x_1,y_2,0,...)$ and $\y_2=(y_2,0,...)$. Similarly $\X'$ is
the set of points $\x_0'=(x_0',x_1',x_2',0,...)$,
$\x_1'=(x_1',x_2',0,...)$, $\x_2'=(x_2',0,...)$,
$\y_0'=(x_0',y_1',y_2',0,...)$, $\y_1'=(y_1',y_2',0,...)$ and
$\y_2'=(y_2',0,...)$.
Hence the systems $(\X,\tal)$  and $(\X',\tal')$ are given by the same diagram \\
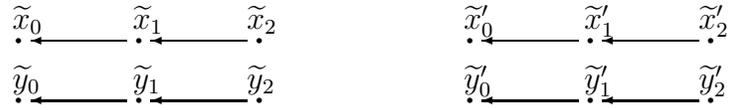
\begin{picture}(150,22)(-40,-2)
\put(15,12){\circle*{1}} \put(14,14){$\x_0$}
\put(33,12){\vector(-1,0){16}}

 \put(15,2){\circle*{1}} \put(14,4){$\y_0$}
\put(33,2){\vector(-1,0){16}}

 \put(35,12){\circle*{1}}\put(34,14){$\x_1$}
\put(53,12){\vector(-1,0){16}}

 \put(35,2){\circle*{1}} \put(34,4){$\y_1$}
\put(53,2){\vector(-1,0){16}}

 \put(55,12){\circle*{1}} \put(53,14){$\x_2$}

 \put(55,2){\circle*{1}} \put(53,4){$\y_2$}

 \put(90,12){\circle*{1}} \put(89,14){$\x_0'$}
\put(108,12){\vector(-1,0){16}}

 \put(90,2){\circle*{1}} \put(89,4){$\y_0'$}
\put(108,2){\vector(-1,0){16}}

 \put(110,12){\circle*{1}}\put(109,14){$\x_1'$}
\put(128,12){\vector(-1,0){16}}

 \put(110,2){\circle*{1}} \put(109,4){$\y_1'$}
\put(128,2){\vector(-1,0){16}}

 \put(130,12){\circle*{1}} \put(128,14){$\x_2'$}

 \put(130,2){\circle*{1}} \put(128,4){$\y_2'$}
\end{picture}
\end{exm}
\subsection{Topological Markov chains, projective limits and hiperbolic attractors}\label{subsection}
 It is not hard to give an example of a  partial   dynamical system which is not a reversible extension of any
'smaller' dynamical system, though  its dynamics is implemented by a partial homeomorphism.
 Yet  many reversible dynamical systems arise from irreversible ones as their reversible extensions.
  In order to see that we recall first  the  topological Markov chains and  then the hyperbolic attractors.\\
  Let $A=(A(i,j))_{i,j\in\{1,...,N\}}$ be  a square matrix with entries in $\{0,1\}$, and
  such that no row  of $A$ is identically zero. We associate with $A$ two dynamical systems $(X_A,\sigma_A)$
   and $(\overline{X}_A,\overline{\sigma}_A)$.
 The \emph{one-sided  Markov subshift}\label{one-sided Markov} $\sigma_A$ acts on the compact space
 $X_A=\{(x_k)_{k\in \N}\in \{1,...,N\}^{\N}:A(x_k,x_{k+1})=1,\,\,k\in\N\}$
  (the topology on $X_A$ is the one inherited from the Cantor space $\{1,...,N\}^{\N}$) by the rule
 $$(\sigma_A x)_k=x_{k+1},\qquad \mathrm{for}\,\, k\in \N,
 \,\mathrm{and}\,\, x\in X_A .$$
Unless $A$ is a permutation matrix $\sigma_A$ is not injective,
and $\sigma_A$ is onto if and only if every column of $A$ is
non-zero. The \emph{two-sided Markov subshift}
$\overline{\sigma}_A$ acts on the compact space
$\overline{X}_A=\{(x_k)_{k\in \Z}\in
\{1,...,N\}^{\Z}:A(x_k,x_{k+1})=1,\,\,k\in\Z\}$ and is defined by
 $$(\overline{\sigma}_A x)_k=x_{k+1},\qquad \mathrm{for}\,\, k\in \Z,
 \,\mathrm{and}\,\, x\in \overline{X}_A .$$
Mapping $\overline{\sigma}_A$ is what is called a
\emph{topological Markov chain} and  abstractly can be defined as
an expansive homeomorphism of a completely disconnected compactum.
\\
If we  assume that not only the  rows  of $A$  but also the columns are not identically zero then $\sigma_A$ is onto and we have

\begin{prop}\label{topologicalmarkovchain}
Let  $A$ have no zero columns and let $(\X_A,\widetilde{\sigma}_A)$ be the reversible extension of  $(X_A,\sigma_A)$.
 Then
$$
(\X_A,\widetilde{\sigma}_A)\cong (\overline{X}_A,
\overline{\sigma}_A).
$$
\end{prop}
\begin{Proof}
Let $\x\in\X_A$. Then $\x=(x_0,x_1,...)$ where
$x_n=(x_{n,k})_{k\in\N}\in X_A$  is such that
$\sigma_A^m(x_{n+m})=x_{n}$, $n,m\in \mathbb{N}$. Thus,  $x_{n+m,k+m}=x_{n,k},$ for $k,\,n,\,m\in\N$, or in other
words
$$
m_1-n_1=m_2-n_2\Longrightarrow x_{n_1,m_1}=x_{n_2,m_2}.$$
 We see  that the sequence $\overline{x}\in \overline{\sigma}_A$
such that $\overline{x}_k=x_{n,m}$, $m-n=k\in \mathbb{Z}$, carries
the full information about $\x$. Hence, defining  $\Upsilon$ by
the formula
\begin{equation}\label{Upsilon}
\Upsilon(\x)=(...,x_{n,0},...,x_{1,0},
\dot{x}_{0,0},x_{0,1},...,x_{0,n},...)
\end{equation}
where dot over $x_{0,0}$ denotes  the zero entry, we get
injective mapping $\Upsilon:\X_A\rightarrow \overline{X}_A$. It is
evident that $\Upsilon$ is
surjective and can be readily checked that it is  also  continuous, whence $\Upsilon$ is a homeomorphism.\\
Finally  let us recall that
$\widetilde{\sigma}_A(\x)=(\sigma_A(x_0),x_0,...)$ and
$\sigma_A(x_0)=(x_{0,1},x_{0,2},...)$ and thus
$$
\Upsilon(\widetilde{\sigma}_A(\x))=(...,x_{n,0},...,x_{1,0},
x_{0,0},\dot{x}_{0,1},...,x_{0,n},...)=\overline{\sigma}_A(\Upsilon(\x))
$$
which says that $(\X_A,\widetilde{\sigma}_A)$ and
$(\overline{\sigma}_A,\overline{X}_A)$ are topologically conjugate
by $\Upsilon$.
\end{Proof}

  The  shift $\sigma_A$ has a clopen image  of the form
$$
\sigma_A(X_A)=\{(x_k)_{k\in \N}\in X_A: \sum_{x=1}^N
A(x,x_0)>0\}.
$$ 
Thus, if $A$   has  at least one  zero column then $\sigma_A$ is  not onto,
 and since $\overline{\sigma}_A$ is always
onto,  systems $(\X_A,\widetilde{\sigma}_A)$ and $(\overline{X}_A,
\overline{\sigma}_A)$ cannot be  conjugate.
In fact, $\X_A=\bigcup_{n\in\N} X_{A,n} \cup X_{A,\infty}$, see  Definition \ref{widetildeXdefn}, 
and in the same manner as in the proof of Proposition
\ref{topologicalmarkovchain} we may define  a homeomorphism
$\Upsilon$ from the subset $X_{A,\infty}=\{(x_0,x_1,...): x_n\in X_A,\,
\sigma_A(x_{n})=x_{n-1},\, \,n\geqslant 1\}$ onto  $\overline{X}_A$ such that 
$\Upsilon(\widetilde{\sigma}_A(\x))=\overline{\sigma}_A(\Upsilon(\x))$ for $\x\in X_{A,\infty}$, that is
$$
(X_{A,\infty},\widetilde{\sigma}_A)\cong (\overline{X}_A,
\overline{\sigma}_A) .
$$
In order to build a homeomorphic image of the whole space $\X_A$ we need to  add countable number
of  components to $\overline{X}_A$. We may do it by embedding $\overline{X}_A$ into another space $\overline{X}_{A'}$
 associated with the larger alphabet $\{0,1,...,N\}$ and a larger matrix $A'=(A'(i,j))_{i,j\in\{0,1,...,N\}}$.
\begin{prop}
Let $A'=(A'(i,j))_{i,j\in\{0,1,...,N\}}$ be given by
$$
A'(i,j)=\begin{cases}
        A(i,j), &  \textrm{if }\,\, i,j\in \{1,...,N\},\\
        1 ,     &  \textrm{if }\,\, i=0 \textrm{ and either }j=0  \textrm{ or  }j\textrm{-th column of } A
         \textrm{ is zero }  ,\\
        0,       & \textrm{ otherwise},
\end{cases}
$$
and let $\overline{X}_{A'}^{+}=\{(x_k)_{k\in \Z}\in\overline{X}_{A'}
:x_0\neq 0\}$. Then  $\overline{\sigma}_{A'}(\overline{X}_{A'}^{+})\subset \overline{X}_{A'}^{+}$ and
$$
(\X_A,\widetilde{\sigma}_A)\cong (\overline{X}_{A'}^+,
\overline{\sigma}_{A'}).
$$
\end{prop}
\begin{Proof}
Let us treat   $\overline{X}_A$  as a subset of
$\{0,1,...,N\}^{\Z}$ and recall that $\X_A=\bigcup_{n\in\N}
X_{A,n} \cup X_{A,\infty}$, and we have   the homeomorphism
$\Upsilon:X_{A,\infty}\rightarrow \overline{X}_A$, cf.
\eqref{Upsilon}.  We put 
$$
\Upsilon(\x)=(...,0,0,x_{n-1,0},...,x_{1,0},
\dot{x}_{0,0},x_{0,1},...,x_{0,n},...)$$
for $\x=(x_0,x_1,...,x_{n-1},0,...)\in X_{A,n}$
where $x_m=(x_{m,k})_{k\in\N}\in X_A$,  $m=0,...,n-1$. 
 In the same fashion as in
the proof of  Proposition \ref{topologicalmarkovchain} one checks
that  the mapping $\Upsilon:\X_A\rightarrow \{0,1,...,N\}^{\Z}$ is
injective and the equality
$\Upsilon(\widetilde{\sigma}_{A}(\x))=\overline{\sigma}_{A'}(\Upsilon(\x))
$ holds for every $\x\in \X_A$. Moreover,  $X_{A,n}$  is  mapped
by $\Upsilon$ onto the set 
$$
 \{(x_k)_{k\in \Z}\in \{0,1,...,N\}^{\Z}:
x_k=0,\,\,k< -n; \sum_{i=1}^N A(i,x_n)=0; \,A(x_k,x_{k+1})=1,\,\,k
> - n\}
$$
denoted by $\overline{X}_{A,n}$. Thus $
(\X_A,\widetilde{\sigma}_A)\cong (\bigcup_{n\in
\N}\overline{X}_{A,n}\cup \overline{X}_A, \overline{\sigma}_{A'})$. It is clear that
$\overline{X}_{A'}^{+}=\bigcup_{n\in
\N}\overline{X}_{A,n}\cup \overline{X}_A$ and hence the proof is complete.
\end{Proof}
 The proof of Proposition \ref{topologicalmarkovchain} is  actually  the proof  of  the probably known
fact that  if $\sigma_A$ is onto,  then $\overline{X}_A$ is the
projective (inverse) limit of the projective sequence
$X_A\stackrel{\sigma_A}{\longleftarrow} X_A
\stackrel{\sigma_A}{\longleftarrow} ...$ . Let us pick out the relationship between
$\al$-extensions and projective limits.
For that purpose (and also for future needs) we
 introduce some  terminology.
\\
As $\tal$ is a partial
homeomorphism, we  denote by $\TDelta_{n}$  the domain of
$\tal^n$,  $n\in \mathbb{Z}$. For $n\in \N$ we have
$$
\TDelta_{n}=\{\x=(x_0,x_1,...)\in \X: x_0 \in
\Delta_{n}\}=\Phi^{-1}(\Delta_n),
$$
$$
\TDelta_{-n}=\{\x=(x_0,x_1,...)\in \X: x_{n} \neq 0\}\subset
\Phi^{-1}(\Delta_{-n}).
$$
 With the help of
 $\Phi$ and $\tal^{-1}$, we define the family of projections
 $\Phi_n=\Phi\circ\tal^{-n}:\TDelta_{-n} \longrightarrow \Delta_n$, for $n\in\N$. We have
$$\label{Phi_n}
\Phi_n(\x)=(\Phi\circ\tal^{-n})(x_0,x_1,...,x_n,...)=x_n,\quad \qquad \x\in \TDelta_{-n}.
$$
Since  $X_\infty=\bigcap_{n=1}^\infty \TDelta_{-n}$, the mappings $\Phi_n$ are well
defined on $X_\infty$, and the following statements are straightforward.
\begin{prop}\label{projectivelimit}
 The system  $(X_\infty,\Phi_n)_{n\in\N}$ is the  projective
limit  of the
projective sequence $(\Delta_n,\al_n)_{n\in\N}$ where
$\al_n=\al|_{\Delta_n}$: $X_\infty =\underleftarrow{\,\,\mathrm{lim}\,}(\Delta_n,\al_n)$.
\end{prop}
\begin{corl} If $\al$ is onto, then $\X =\underleftarrow{\,\,\mathrm{lim}\,}(\Delta_n,\al_n)$.
\end{corl}
The above statement provides us with  many interesting examples of reversible extensions because    projective (inverse) limit spaces commonly appear as attractors in dynamical systems (this was observed for the first time by  R. F. Willams,   see \cite{Williams}). We recall here the  classic example.  
\begin{exm}[Solenoid]\label{solenoid} 
Let $S^1=\{z\in\C: |z|=1\}$ be the unit circle in the complex plane and let $\al$ be the expanding endomorphism
of $S^1$ given by
 $$\al(z)=z^2, \qquad z\in S^1 .$$ 
  Then the projective limit $\underleftarrow{\,\,\lim\,}(S^1,\al)$ is homeomorphic to the Smale's  solenoid,
 that is an attractor of the mapping $F$ acting on the solid torus $\TT=S^1\times D^2$, where $D^2=\{z\in \C: |z|\leq 1\}$, by $
F(z_1,z_2)=(z_1^2,\lambda z_2+ \frac{1}{2}\,z_1)$
 where $0<\lambda<\frac{1}{2}$ is fixed.
 \begin{center}\setlength{\unitlength}{0.7mm}
\begin{picture}(110,70)(15,0)\thicklines
%\put(0,2){\includegraphics[angle=0, scale=1.3]{solenoid.eps}} 
%Donuthole 
\qbezier(29,45)(43,40)(59,45)
\qbezier(62,50)(64,47.5)(59,45)
\qbezier(29,45)(24,47.5)(26,50)
\qbezier(27.5,45.8)(34,52)(44,52)
\qbezier(44,52)(55,52)(60,45.5)
%Donutcontour
\qbezier(24,6.3)(43,0)(65,6.7)
\qbezier(65,6.7)(68,7.5)(73,10.5)
\qbezier(24,6.3)(21,7.1)(15.5,10)
\qbezier(85.2,28)(83.5,18)(73,10.5)
\qbezier(3.1,28)(4.5,18)(15.5,10)
\qbezier(85.2,28)(88.5,43)(79.2,54)
\qbezier(3.1,28)(-0.5,43)(9,54)
\qbezier(44,66.3)(66,66.5)(79.2,54)
\qbezier(44,66.3)(21,66.5)(9,54)
%insideDonut
\thinlines
\qbezier(19.5,33.5)(16,32.7)(15,31)
\qbezier(39,27.6)(21,26)(14.2,31.7)
\qbezier(39,27.6)(52.5,27.6)(68,37.5)
\qbezier(71,43)(72,41)(68,37.5)
\qbezier(71,43)(68,48)(59,49.5)
\qbezier(43,50)(50,51)(59,49.5)
\qbezier(43,50)(32,48.5)(28,43)
\qbezier(40,28)(21,35)(28,43)
\qbezier(23.2,28.1)(18.5,33)(19,40)
\qbezier(33.2,57)(19,50)(19,40)
\qbezier(33.2,57)(47,62)(59,59)
\qbezier(76.2,50)(69,58)(59,59)
\qbezier(76.2,50)(79,44)(78.5,40)
\qbezier(62,23.5)(76,29)(78.5,40)
\qbezier(62,23.5)(47,18)(34,18)
\qbezier(8,25.5)(14,18)(34,18)
\qbezier(8,25.5)(4,32)(7,36)
\qbezier(19,43)(11,41)(7,36)
\qbezier(30.4,45)(44,46)(58,42)
\qbezier(67.3,37.3)(64,40)(58,42)
\qbezier(28,35)(44,37)(58,32)
\qbezier(73,30)(76,20)(60,16.3)
\qbezier(40,18)(46,15.6)(60,16.3)
%
%Podpis
\put(108, 44){\textbf{Figure 2.}}

 \put(108, 36){The image of the solid torus}

\put(108, 30){under $F$ is a solid torus which}

 \put(108, 24){wraps twice around itself.}
\end{picture}
\end{center} 
Namely the solenoid is the set $\Sol=\bigcap_{k\in\N}F^k(\TT)$, see Fig. 2, and 
 the reversible extension of $(S^1,\al)$
  is equivalent to $(\Sol,F|_{\Sol})$, see e.g. \cite{introdynsys}.
  
\end{exm}
\subsection{Decomposition of sets in $\X$}
Before the end of this section we introduce  a certain  idea which
 enables us to 'decompose' a subset $\U\subset \X$ into the family $\{U_n\}_{n\in\N}$ of
subsets of $X$. We shall need this device in Section 5.
\begin{defn}\label{U_n}
Let $\U\subset\X$ be a subset of $\al$-extension of $X$ and let
$n\in \N$. We shall call the set
$$
 U_n=\Phi_n(\U\cap \TDelta_{-n})
 $$
 an  \emph{$n$-section} of $\U$.
\end{defn}
\noindent It is evident  that if $\{U_n\}_{n\in\N}$ is the family
of sections of $\U$ then $ \U$ is a subset of $\Big(U_0\times (U_1
\cup\{0\})\times ...\times (U_n \cup\{0\})\times ...\Big)\cap\X $
but  in general
 the opposite  relation does not hold (see Example \ref{exm x_1-x_0}). Fortunately
 we have the following statements true.
\begin{prop}\label{Implication}
If $\al$ is injective on the inverse image of
$\Delta_{-\infty}:=\bigcap_{n\in\N}\Delta_{-n}$, that is for every
point $x\in \Delta_{-\infty}$ we have $|\al^{-1}(x)|=1$, then  for
every subset $\U\subset\X$   we have
\begin{equation}\label{widetildeU1} \U=\Big(U_0\times (U_1
\cup\{0\})\times ...\times (U_n \cup\{0\})\times ...\Big)\cap\X
\end{equation}
where $U_n$ is the $n$-section of $\U$, $n\in \N$.
\end{prop}
\begin{Proof}
Let $\x=(x_0,x_1,...)\in (U_0\times (U_1 \cup\{0\})\times
...\times (U_n \cup\{0\})\times ...)\cap\X$. We show that $\x\in
\U$. Indeed, if $\x\notin X_\infty=\bigcap_{n\in\N}\TDelta_{-n}$,
then  $\x=(x_0,...,x_N,0,...)$ where $x_N\in \Delta_N\setminus
\Delta_{-1}$, and thus $\x$ is uniquely determined by $x_N$. As $x_N\in U_N=\Phi_N(\U\cap
\TDelta_{-N})$,  $\U$ must contain $\x$. \\
If
$\x\in X_{\infty}$ then  $x_n\in\Delta_{-\infty}$ for all $n\in \N$, and $\x$ is
uniquely determined by $x_0\in U_0$. Thus $\x\in \U$.
\end{Proof}

\begin{thm}\label{closed set brackets}
Let $(X,\al)$ be a partial dynamical system such that $\Delta_{-1}$ is open. Then
 every closed subset $\widetilde{V}\subset \X$ is uniquely
determined by its sections $\{V_n\}_{n\in\N}$ via formula
(\ref{widetildeU1}).
\end{thm}
\begin{Proof}
 Let $\widetilde{V}\subset \X$ be closed, that is compact (see Theorem \ref{compactnes}),  and
let  $\x=(x_0,x_1,...)\in (V_0\times (V_1 \cup\{0\})\times
...\times (V_n \cup\{0\})\times ...)\cap\X$. We show that $\x\in \V$. If
$\x\notin X_\infty$ then (see the argument in the proof
of Proposition \ref{Implication})  we immediately get
$\x\in\widetilde{V}$. Thus we only need to consider the case
when
$\x\in X_\infty$.
For that purpose we  define
$\widetilde{D}_n=\{\widetilde{y}=(y_0,y_1,...)\in \X:
y_n=x_n\}\cap\widetilde{V}$,  $n\in\N$. Clearly
$\{\widetilde{D}_n\}_{n\in\N}$ is the decreasing family of closed
nonempty subsets of the compact set $\widetilde{V}$. Hence
$$ \bigcap_{n\in\N}\widetilde{D}_n=\{\x\}\in \widetilde{V}$$
and the proof is complete.
\end{Proof}
\begin{exm}\label{exm x_1-x_0}
For the sake of  illustration of the preceding statements and to see that they cannot
 be sharpen let us consider a dynamical system given by the diagram\\
\setlength{\unitlength}{.8mm}
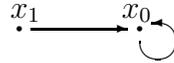
\begin{picture}(35,13)(-85,0)
\put(10,5){\circle*{1}} \put(8.5,7){$x_1$} \put(30,5){\circle*{1}}
\put(27,7){$x_0$}
 \put(33,3){\oval(
6,6)[r]} \put(33,3){\oval(6,6)[bl]} \put(33,6){\vector(-1,0){1}}
\put(12,5){\vector(1,0){16}}
\end{picture}\\
or equivalently by relations $X=\Delta_1=\{x_0,x_1\}$,
$\al(x_1)=\al(x_0)=x_0$ (or equivalently let $\al=\sigma_A$ and
$X=X_A$ where $A=\left(\begin{array}{c c}
 1   &     0  \\
 1   &     0  \\
\end{array}\right)$). Then $\Delta_{-\infty}=\{x_0\}$ and $|\al^{-1}(x_0)|=2$, therefore Proposition
\ref{Implication} can not be used. The space $\X$ consists of
elements $\x_n=(\underbrace{x_0,...,x_0,x_1}_n,0,...)$, $n\in\N$,
and $\x_\infty=(x_0,...,x_0,...)$. Hence it is convenient to
identify $\X$ with the compactification
$\overline{\N}=\N\cup\{\infty\}$
 of the discrete space $\N$. Under this identification  $\tal$ is given by
$$
\tal(n)=n+1,\,\, n\in \N,\qquad\qquad \tal(\infty)=\infty.
$$
It is clear that  all the  sections of the subset  $\N\subset \overline{\N}$ are  equal to $X$.
  As $\N$  is not closed in $\overline{\N}$,  Theorem \ref{closed set brackets} does not work here. Indeed, we have
 $$
\overline{\N}=\Big(X\times (X \cup\{0\})\times ...\times (X
\cup\{0\})\times ...\Big)\neq \N.
 $$ \end{exm}

\section{Covariant representations and their coefficient algebra}
The aim of this section is to study the interrelations between the
covariant representations of $C^*$-dynamical systems corresponding
to $(X,\al)$ and its reversible extension $(\X,\tal)$. Our main tool will be Theorem \ref{ideals2} and hence, cf.
  Proposition \ref{prop(O)sition},
\\

\emph{from now on we  shall
 always  assume that the image $\Delta_{-1}$ of
the partial mapping $\al$ is open.}
\\

\noindent First we show that the algebra $\B=C(\X)$ possesses a certain universal property with
respect to  covariant faithful representations of $(X,\al)$. Afterwards, we construct a dense
$*$-subalgebra of $\B$ with the help of which we investigate the
structure of $\B$ and endomorphisms induced by $\tal$ and
$\tal^{-1}$. Finally we show that there is a one-to-one
correspondence between the covariant faithful representations of
$(X,\al)$ and $(\X,\tal)$, and in the case $\al$ is onto this
correspondence is true for all covariant (not necessarily
faithful) representations.
\subsection{Definition and basic result}
 Let us recall that $\A=C(X)$ and $\delta$ is combined
with $\al$ by (\ref{e3.0}). We  denote by $C_{K}(X)$  the algebra
of continuous functions on $X$ vanishing outside a set $K\subset
X$.
 We start with the definition of covariant representation,
 cf. \cite{Pedersen}, \cite{exel}, \cite{mcclanachan}, \cite{Adji_Laca_Nilsen_Raeburn}.
\begin{defn}\label{cov-rep-def}
 A \emph{covariant
representation} of a $C^*$-dynamical system $(\A,\delta)$, or of
the  partial dynamical system $(X,\al)$, is a triple $(\pi,U,H)$
where $\pi:\A\rightarrow L(H)$ is a representation of $\A$ on
Hilbert space $H$
 and $U\in L(H)$ is  a partial isometry whose initial space is $\pi(C_{\Delta_{-1}}(X))H$ and whose final space is
 $\pi(C_{\Delta_{1}}(X))H$. In addition it is required that
$$
U\pi(a)U^*=\pi(\delta(a)),\qquad \,\,\mathrm{for}\, \,a\in \A.
$$
If the representation $\pi$ is faithful we call the  triple
$(\pi,U,H)$ a \emph{covariant faithful representation}. Let
CovRep$(\A,\delta)$ denote the set of all covariant representations
and CovFaithRep$(\A,\delta)$  the set of all covariant faithful
representations of $(\A,\delta)$.
\end{defn}
\begin{rem}\label{covreprem}
As  $\Delta_1$ and $\Delta_{-1}$ are clopen, the projections
$P_{\Delta_{1}}$ and  $P_{\Delta_{-1}}$ corresponding to the
characteristic functions $\chi_{\Delta_{1}}$ and
$\chi_{\Delta_{-1}}$ belong to $\A$. Thus for every covariant
representation $(\pi,U,H)$ we see that $UU^*=\pi(P_{\Delta_{1}})$
and $U^*U=\pi(P_{\Delta_{-1}})$ belong to $\pi(\A)$, cf.
(\ref{condition}).
\end{rem}
\noindent We shall see in Corollary \ref{gupiecorl} that  every
$C^*$-dynamical system  $(\A,\delta)$ possesses a
covariant faithful representation, and hence the sets $\mathrm{CovRep}\,(A,\delta)$ and
$\mathrm{CovFaithRep}\,(\A,\delta)$ are non-empty. \\
Now we  reformulate the main result of
\cite{maxid}  in terms of covariant representations. The point is
that for every   covariant  representation $(\pi,U,H)$ of
$(\A,\delta)$ the condition \eqref{condition} holds, whence  if $\pi$ is
faithful then  by Theorem \ref{ideals2} the maximal ideal
space of the $C^*$-algebra $C^*(\bigcup_{n=0}^\infty U^{*n} \pi(\A) U^n)$ is homeomorphic to
$\al$-extension $\X$ of $X$.
 \begin{thm}\label{coefalgebra}
 Let $(\pi,U,H) \in \mathrm{CovFaithRep}\,(\A,\delta)$ and let $\X$ be the $\al$-extension of $X$.
Then
$$
C^*\big(\bigcup_{n=0}^\infty U^{*n} \pi(\A) U^n\big)\cong C(\X).
$$
In other words, the  coefficient $C^*$-algebra of $C^*(\pi(\A),U)$ is
isomorphic to the algebra of continuous functions on
$\alpha$-extension of $X$, cf. \cite{lebiedodzij}.
Moreover this isomorphism   maps an operator of the form
$\pi(a_0)+U^*\pi(a_1)U+...+U^{*N}\pi(a_N)U^N$, where $a_0,a_1,...,a_N\in \A$,  onto
a function $b\in C(\X)$ such that
$$
 b(\widetilde{x})=a_0(x_0)+a_1(x_1)+...+a_N(x_N),
$$
where $ \widetilde{x}=(x_0,...)\in \X$ and we set $a_n(x_n)=0$
whenever $x_n=0$.
\end{thm}
\begin{Proof}
By  Theorem \ref{ideals2} the maximal ideal space of
$C^*(\bigcup_{n=0}^\infty U^{*n} \pi(\A) U^n)$ is homeomorphic to $\X$. Hence $
C^*(\bigcup_{n=0}^\infty U^{*n} \pi(\A) U^n)\cong C(\X) $. Taking into account formulas
(\ref{xi2}),(\ref{xi000'}) we obtain the postulated form of this
isomorphism. Indeed, if $\x=(x_0,x_1,...)$ is a character on $C^*(\bigcup_{n=0}^\infty U^{*n} \pi(\A) U^n)$ then
$$
\x(\sum_{n=0}^NU^{*n}\pi(a_n)U^n)=\sum_{n=0}^N \x(U^{*n}\pi(a_n)U^n)=\sum_{n=0}^N \xi_{\x}^n(a_n)=\sum_{n=0}^N a_n(x_n)
$$
where  $\xi^n_{\widetilde{x}}(a)=\x(U^{*n}\pi(a)U^n)$, cf. \eqref{xi2}, and $x_n=0$  whenever $\xi_{\x}^n\equiv 0$.
\end{Proof}
 From the above it follows that $\B=C(\X)$ can be regarded as the universal (in fact unique)
  coefficient $C^{*}$-algebra for covariant faithful
  representations. In case $\delta$ is injective (that is $\al$ is onto), the universality of
  $\B$ is  much 'stronger', see Proposition \ref{directlimit}.
\subsection{Coefficient $^*$-algebra}

  We shall  present now a certain dense $^*$-subalgebra of $\B=C(\X)$, a coefficient $^*$-algebra  which
    frequently  might be more convenient to work with.
   The plan is to construct an algebra $B_0\subset
   l^1(\mathbb{N},\A)$ and
    then take the quotient of it by  certain ideal. The result will be naturally isomorphic to
    a  $^*$-subalgebra of $\B$.\\
  First, let us observe that if we set
$\A_n:=\delta^n(1)\A$, $n\in \N$, then we obtain a decreasing
family $\{\A_n\}_{n\in\N}$, of  closed two-sided ideals in $\A$. Since
the operator $\delta^n(1)$ corresponds to the characteristic function $\chi_{\Delta_n}\in C(X)$,
one can consider  $\A_n$ as $C_{\Delta_n}(X)$. Let $
B_0$  denote the set consisting of sequences
$a=\{a_n\}_{n\in \N}$ where $a_n\in \A_n$, $n\in \N$,
 and only  a finite number of functions $a_n$  is non zero. Namely
$$
B_0=\{a\in \prod_{n=0}^\infty \A_n: \exists_{N>0}\,
\forall_{n>N}\quad a_n\equiv 0\}.
$$
Let $a=\{a_n\}_{n\geqslant 0}$, $b=\{b_n\}_{n\geqslant 0}\in
B_0$ and $\lambda\in\mathbb{C}$. We define the addition,
multiplication by scalar, convolution multiplication and
involution on $ B_0$  as follows
\begin{equation}\label{add} (a+b)_n=a_n+b_n, \end{equation}
\begin{equation}\label{mulscal}
 (\lambda a)_n=\lambda a_n,
\end{equation} \begin{equation}\label{mul}
 (a\cdot b)_n=a_n\sum_{j=0}^n \delta^j(b_{n-j})+b_n \sum_{j=1}^n
\delta^j(a_{n-j}), \end{equation} \begin{equation}\label{invol}
 (a^*)_n=\overline{a}_n.
\end{equation}
 These operations are well defined
and seems very familiar, except maybe the multiplication of two elements
from $B_0$. We point out here that the index in one
of the sums of  (\ref{mul}) starts running from $0$.
\begin{prop}\label{algebrawithinvolution}
The set $B_0$ with operations (\ref{add}),
(\ref{mulscal}), (\ref{mul}), (\ref{invol}) becomes a commutative
  algebra with involution.
\end{prop}

\begin{Proof}
It is clear that operations (\ref{add}), (\ref{mulscal}) define
the structure of vector space on $B_0$ and that
operation (\ref{invol}) is an involution. The rule (\ref{mul}) is
less easy to show up its properties. Commutativity and
distributivity can be checked easily but  in order to
prove the associativity we must strain ourselves quite a lot.\\
Let $a$, $b$, $c \in B_0$. Then
$$
((a\cdot b)\cdot c)_n=(a\cdot b)_n \sum_{j=0}^n
\delta^j(c_{n-j})+c_n \sum_{j=1}^n \delta^j((a\cdot b)_{n-j}),
$$
$$
(a\cdot (b\cdot c))_n=a_n \sum_{j=0}^n \delta^j((b\cdot c)_{n-j})+
(b\cdot c)_n \sum_{j=1}^n \delta^j(a_{n-j}),
$$
where
$$
(a\cdot b)_n \sum_{j=0}^n \delta^j(c_{n-j})= [a_n\sum_{k=0}^n
\delta^k(b_{n-k})+b_n \sum_{k=1}^n
\delta^k(a_{n-k})]\sum_{j=0}^n\delta^j(c_{n-j})
$$
$$
=a_n\sum_{k,j=0}^n \delta^k(b_{n-k})\delta^j(c_{n-j})+b_n \sum_{
k=0, j=1 }^n \delta^k(a_{n-k})\delta^j(c_{n-j})
$$
and
$$
c_n \sum_{j=1}^n \delta^j((a\cdot b)_{n-j})=c_n
\sum_{j=1}^n\delta^j(a_{n-j}\sum_{k=0}^{n-j}
\delta^k(b_{n-j-k})+b_{n-j} \sum_{k=1}^{n-j} \delta^k(a_{n-j-k}))
$$
$$
=c_n \sum_{j=1}^n\delta^j(a_{n-j})\sum_{k=j}^{n}
\delta^k(b_{n-k})+\delta^j(b_{n-j})\sum_{k=j+1}^{n}
\delta^k(a_{n-k})) =c_n\sum_{ k=1, j=1 }^n
\delta^j(a_{n-j})\delta^k(b_{n-k}).
$$
Simultaneously  by analogous computation
$$
a_n \sum_{j=0}^n \delta^j((b\cdot c)_{n-j})=a_n\sum_{k,j=0}^n
\delta^k(b_{n-k})\delta^j(c_{n-j})
$$
$$
(b\cdot c)_n \sum_{j=1}^n \delta^j(a_{n-j})=b_n \sum_{ k=0, j=1
}^n \delta^k(a_{n-k})\delta^j(c_{n-j})+c_n\sum_{ k=1, j=1 }^n
\delta^j(a_{n-j})\delta^k(b_{n-k}).
$$
Thus,  $((a\cdot b)\cdot c)_n=(a\cdot (b\cdot c))_n$ and the
$n$-th entry of the sequence $a\cdot b\cdot c$ is of the form
$$
a_n\sum_{k,j=0}^n \delta^k(b_{n-k})\delta^j(c_{n-j})+b_n \sum_{
k=0, j=1 }^n \delta^k(a_{n-k})\delta^j(c_{n-j})+c_n\sum_{ k=1, j=1
}^n \delta^k(a_{n-k})\delta^j(b_{n-j}).
$$
\end{Proof}
Let us define a morphism
$\varphi:B_0\rightarrow \B$. To this end, let
$a=\{a_n\}_{n\in \N}\in B_0$ and
$\x=(x_0,x_1,...)\in \X$. We set
\begin{equation}\label{phimap}
\varphi(a)(\x)=\sum_{n=0}^\infty a_n(x_n),
\end{equation}
where  $a_n(x_n)=0$ whenever $x_n=0$. The mapping $\varphi$ is
well defined as only a finite number of functions $a_n$, $n\in
\N$, is non zero.
\begin{thm}
The mapping $\varphi:B_0\rightarrow \B$ given by
(\ref{phimap}) is a morphism of algebras with involution and the
image of $\varphi$ is dense in $\B$, that is
$$\overline{\varphi(B_0)}=\B.$$
\end{thm}
\begin{Proof} It is clear that $\varphi$ is a linear mapping preserving an involution.
We show that $\varphi$ is multiplicative.
 Let $a,$ $b\in
B_0$ and $\x=(x_0,x_1,...)\in \X$ and let $N>0$ be
such that for every $m > N$ we have $a_m=b_m=0$. Using the fact
that $\al^j(x_n)=x_{n-j}$ we obtain
$$
 \varphi(a\cdot b)(\x)=\sum_{n=0}^{N}(a\cdot b)_n(x_n)
 =\sum_{n=0}^{N}\Big[a_n\sum_{j=0}^n \delta^j(b_{n-j})+b_n
 \sum_{j=1}^n
\delta^j(a_{n-j})\Big](x_n)
$$
$$
=\sum_{n=0}^{N}\big[a_n(x_n)\sum_{j=0}^n b_{n-j}(x_{n-j})+b_n(x_n)
 \sum_{j=1}^{n}
a_{n-j}(x_{n-j})\big]
$$
$$
=\sum_{n=0,j=0}^{N} a_n(x_n) b_{j}(x_j)=\sum_{n=0}^{N}a_n(x_n)
\cdot \sum_{j=0}^{N}
 b_{j}(x_j)= \big[\varphi(a)\cdot\varphi(b)\big](\x).
$$
To prove that  $\varphi(B_0)$ is dense in $\B=C(\X)$
we use the Stone-Weierstrass theorem.
  It is clear that $\varphi(B_0)$ is a self-adjoint subalgebra of
 $\B$ and as  $a=(1,0,0,...)\in
 B_0$, we get  $\varphi(a)=1\in
 \B$, that is $\varphi(B_0)$ contains the identity.
 Thus, what we  only  need   to prove is that
 $B_0$  separates points of $\X$.\\
 Let  $\x=(x_0,x_1,...)$ and
 $\widetilde{y}=(y_0,y_1,...)$ be two distinct points of $\X$.
Then  there exists $n\in\mathbb{N}$ such that   $x_n\neq
 y_n$   and  by Urysohn's lemma there  exists a function
 $a_n\in C_{\Delta_n}(X\cup\{0\})$ such that $a_n(x_n)=1$ and
 $a_n(y_n)=0$. Taking  $a\in B_0$ of
 the form
 $$a=(\underbrace{0,...,0}_{n},a_n,0,...)$$
we see that $\varphi(a)(\x)=1$ and $ \varphi(a)(\widetilde{y})=0$.
Thus the proof is complete.
 \end{Proof}
Let us consider the quotient space $B_0/
\mathrm{Ker}\,\varphi$ and the quotient mapping
$\phi:B_0/ \mathrm{Ker}\,\varphi\rightarrow \B_0$,
that is $\phi(a+\mathrm{Ker}\,\varphi)=\varphi(a)$. Clearly $\phi$
is an injective mapping onto a dense $^*$-subalgebra of $\B$. In what follows we make use of  the following notation
$$
\B_0:=\phi(B_0/ \mathrm{Ker}\,\varphi )
$$
$$
 [a]:=\phi(a+\mathrm{Ker}\,\varphi), \qquad a\in B_0
$$
\begin{defn}\label{idenitify-defn}
We shall call $\B_0$ the \emph{coefficient}  \emph{$^*$-algebra}
of a dynamical system  $(\A,\delta)$. We will write
$[a]=[a_0,a_1,...]\in \B_0$ for $a=(a_0,a_1,a_2,...)\in
B_0$.
\end{defn}
 The natural injection $ \A\ni a_0\longrightarrow
[a_0,0,0,...]\in \B_0 $ enables us to treat $\A$ as an
$C^*$-subalgebra of $\B_0$ and hence also of $\B$:
$$
\A\subset \B_0\subset \B,\qquad \overline{\B}_0=\B.
$$
Using the mappings $\Phi_n:\TDelta_{-n} \longrightarrow \Delta_n$ (see
subsection \ref{subsection})  one can  embed into $\B_0$ not only $\A$
but all the subalgebras $\A_n=C_{\Delta_{n}}(X)$, $n\in\N$.
Indeed, if we define $\Phi_{*n}:\A_n\rightarrow
\B$ to act as follows
$$
\Phi_{*n}(a)=[\underbrace{0,...,0}_n,a,0,...]=\left\{\begin{array}{ll}
a\circ\Phi_n& ,\x\in
\TDelta_{-n}\\
0 & ,\x\notin \TDelta_{-n} \end{array}\right. ,
$$
then clearly  $\Phi_{*n}$ are  injective. Moreover  we have
$C^*(\bigcup_{n\in\N} \Phi_{*n}(\A_n))=\B$, and in
the case $\delta$ is a monomorphism, that is $\al$ is surjective,
$\{\Phi_{*n}(\A_n)\}_{n\in\N}$ forms an increasing family of
algebras and $\B_0=\bigcup_{n\in\N} \Phi_{*n}(\A_n)$. We are
exploiting this fact in the coming proposition.
\begin{prop}\label{directlimit}
If $\delta$ is injective then  $\B$ is  the
direct limit $\underrightarrow{\lim\,\,\,}\A_n $ of the sequence
$(\A_n,\delta_{n})_{n=0}^\infty$  where $\delta_n=
\delta|_{\A_n}$, $n\in\N$.
\end{prop}
\begin{Proof} Let  $\B=\underrightarrow{\lim\,\,\,}\A_n$ be the  direct limit of the sequence
$(\A_n,\delta_{n})_{n=0}^\infty$, and let
$\psi^n:\A_n\rightarrow\B$ be the natural homomorphisms, see for
instance \cite{Murphybook}. It is straightforward to see that the
diagram
$$
\begin{array}{ccc}
 \A_n & \,\,\,\stackrel{\delta_n}{\longrightarrow} &\A_{n+1}  \\

 & \,\,\, \searrow  \Phi_{*n}  &\,\,\,\,\,\,\,\,\, \downarrow \Phi_{*n+1}\\

 &       &      \B
\end{array}
$$
commutes and hence  there exists a  unique homomorphism
$\psi:\B\rightarrow \B$ such that the diagram
$$
\begin{array}{ccc}
 \A_n & \,\,\,\,\,\stackrel{\psi^{n}}{\longrightarrow} &\B  \\

 &  \searrow  \Phi_{*n}  &\,\,\,\, \downarrow \psi\\

 &       &      \B
\end{array}
$$
commutes. It is evident that $\psi$ is a surjection  (
$\bigcup_{n\in\N}\Phi_{*n}(\A_n)$ generates $\B$)
and as $ \Phi_{*n}(\A_n)$ are increasing  and $\Phi$ is injective,
$\psi$ is injective. Therefore $\psi$ is an isomorphism and the
proof is complete.
\end{Proof}
The preceding proposition points out the relationship between our
approach and the approach presented (among the others) by G. J. Murphy in
\cite{Murphy}. We shall discuss this relationship in the sequel, see Remark \ref{murphyrem}.
\\
 Let us now proceed and consider  endomorphisms of $\B=C(\X)$ given by the formulae
\begin{equation} \label{tildedelta}
 \Tdelta(a)(x)=\left\{\begin{array}{ll}
a(\tal(\x))& ,\x\in
\TDelta_1\\
0 & ,\x\notin \TDelta_1 \end{array}\right.
\qquad
 \Tdelta_*(a)(x)=\left\{\begin{array}{ll}
a(\tal^{-1}(\x))& ,\x\in
\TDelta_{-1}\\
0 & ,\x\notin \TDelta_{-1}
\end{array}\right.
\end{equation} where $\tal:\TDelta_1\rightarrow\TDelta_{-1}$ is a canonical
partial homeomorphism
 of $\X$  defined by 
formula (\ref{tildealpha}). What is  important is that
characteristic functions of $\TDelta_1$ and $\TDelta_{-1}$  belong
to $\A\subset \B_0$. Indeed, we have
$ \chi_{\TDelta_{1}}=[\chi_{\Delta_1},0,0,...]$, $\chi_{\TDelta_{1}}=[\chi_{\Delta_{-1}},0,0,...]\in \A$
 (see remark after
Definition \ref{idenitify-defn}).
Furthermore,  the domain $\TDelta_n$ of the mapping $\tal^{n}$,  $n\in
\mathbb{Z}$, is clopen and it is just an easy exercise to check that for $n\in \mathbb{N}$ we have
$
 \chi_{\TDelta_{n}}=[\chi_{\Delta_n},0,0,...]$, $\chi_{\TDelta_{-n}}=[\underbrace{0,0,...,0}_{n},\Delta_{n},...]\in\B_0.
$
In particular $
\chi_{\TDelta_{-1}}=[\chi_{\Delta_{-1}},0,0,...]=[0,\chi_{\Delta_1},0,...]$.
We are now ready to give an `algebraic' description of $\Tdelta$
and $\Tdelta_* $.
\begin{prop}\label{propendomor}
Endomorphisms $\Tdelta$ and  $\Tdelta_* $ preserve $^*$-subalgebra
$\B_0\subset \B$ and  for $a=(a_0,a_1,a_2,...)\in
B_0$ we have
\begin{equation}\label{algdelta}
\Tdelta([a])=[\delta(a_0) + a_1,a_2,a_3,...],
\qquad  \Tdelta_*([a])=[0,a_0\delta(1),a_1\delta^2(1),...].
 \end{equation}
\end{prop}
\begin{Proof}
 Let $\x=(x_0,x_1,x_2,...)\in \X$.
 In order to prove the first equality in (\ref{algdelta}) it is enough to notice  that for $\x\in\TDelta_1$ we have
  $$ \delta([a])(\x)=[a](\tal(\x))=[a](\al(x_0),x_0,x_1,...)=
 a_0(\al(x_0)) + a_1(x_0) + a_2(x_1)+ ...
 $$
$$
=[\delta(a_0) + a_1,a_2,a_3,...](x_0,x_1,x_2,...)=[\delta(a_0) +
a_1,a_2,a_3,...](\x)
$$
and for  $\x\notin\TDelta_1$  both sides of the left hand side equation
 in (\ref{algdelta}) are equal to zero.\\
In the same manner we show the  validity of the remaining equality. If
$\x\in\TDelta_{-1}$ then
  $$ \delta_*([a])(\x)=[a](\tal^{-1}(\x))=[a](x_1,x_2,x_3...)=
 a_0(x_1) + a_1(x_2) + a_2(x_3)+ ...
 $$
$$
=[0 ,a_0\delta(1),
a_1\delta^2(1),...](x_0,x_1,x_2,...)=[0,a_0\delta(1),a_1\delta^2(1),...](\x)
$$
and if  $\x\notin\TDelta_1$ then both sides of the right hand side equation
in (\ref{algdelta}) are equal to zero.
\end{Proof}
\begin{exm}\label{exm x_1-x_0-coefalg}
The dynamical system $(X,\al)$ from Example \ref{exm x_1-x_0}
corresponds to the $C^*$-dynamical system $(\A,\delta)$ where
$\A=C(\{x_0,x_1\})$ and $\delta(a)\equiv a(x_0)$. We identify $\X$
with $\overline{\N}$ as we did before. Then, since
$a=[a_0,a_1,...,a_N,0,...]$, $a_k\in\A$, $k=0,...N$, is a continuous
function on $\X=\overline{\N}$ we can regard it as a sequence
which has a limit.
 One readily checks that this sequence has the following form: $a(n)=\sum_{k=0}^{n-1}a_k(x_0)+a_n(x_1)$
  for $n=0,...N$, and $a(n)=\sum_{k=0}^{N}a_k(x_0)$ for $n>N$.
Hence $\B_0$ is the  $^*$-algebra of the eventually  constant sequences,  in particular $\A$
consist of sequences of the form
$(a,b,b,b,...)$, $a,b\in\C$. We have
$$
\B_0=\{ (a(n))_{n\in\N}:\exists_{N\in\N}\,\,
\forall_{n,m>N}\,\,a(n)=a(m)\},\qquad \B=\{ (a(n))_{n\in\N}:\exists_{a(\infty)\in \C }\,\,
\lim_{n\to \infty}a(n)=a(\infty)\}
$$
and  within these identifications $\Tdelta$ is the forward  and $\Tdelta_*$ is the backward shift.
\end{exm}
\subsection{The interrelations  between  covariant representations}
The construction of the $^*$-algebra $\B_0$ enables us to
excavate  the inverse of the  isomorphism from Theorem
\ref{coefalgebra}, and  what is more important it enables  us to realize  that
every covariant faithful representation of  $(\A,\delta)$ gives
rise to a
 covariant faithful representation of
$(\B,\Tdelta)$.
\begin{thm}\label{coefalgebra2}
Let $(\pi,U,H)\in\mathrm{CovFaithRep}\,(\A,\delta)$. Then there exists an
extension $\overline{\pi}$ of $\pi$ onto the
coefficient algebra $\B$ such that
$\overline{\pi}:\B\rightarrow
C^*\big(\bigcup_{n=0}^\infty U^{*n} \pi(\A) U^n\big)$ is an isomorphism defined by
$$
 [a_0,...,a_N,0,...] \longrightarrow
\pi(a_0)+U^*\pi(a_1)U+...+U^{*N}\pi(a_N)U^N.
$$
Moreover,  we  have $(\overline{\pi},U,H)\in
\mathrm{CovFaithRep}\,(\B,\Tdelta)$, that is
\begin{equation}\label{takiesobieequation}
\overline{\pi}(\Tdelta(a))=U\overline{\pi}(a)U^*,\,\qquad
\,\,\overline{\pi}(\Tdelta_*(a))=U^*\overline{\pi}(a)U,\qquad\quad
a\in \B.
\end{equation}
\end{thm}
\begin{Proof}
In view of Theorem \ref{coefalgebra} it is immediate  that
$\overline{\pi}$ is an isomorphism. By Theorem \ref{tu3} and by the
form  of endomorphisms
$\Tdelta$ and $\Tdelta_*$, see (\ref{tildedelta}), we get (\ref{takiesobieequation}).
\end{Proof}
We can give a statement somewhat inverse to the above.
\begin{thm}\label{dobrethm}
Let $(\overline{\pi},U,H)\in \mathrm{CovRep}\,(\B,\Tdelta)$  and let $\pi$ be the restriction
  of $\overline{\pi}$ onto  $\A$. Then $(\pi,U,H)\in\mathrm{CovRep}\,(\A,\delta)$.
 Moreover if $(\overline{\pi},U,H)$ is in $\mathrm{CovFaithRep}\,(\B,\Tdelta)$ then
  $(\pi,U,H)$ is in $\mathrm{CovFaithRep}\,(\A,\delta)$ and the extension of $\pi$ mentioned in
  Theorem \ref{coefalgebra2} coincides with $\overline{\pi}$.
\end{thm}
\begin{Proof}
Recall that for $a\in\A$ we write $[a,0,0,...]\in
\B$ and thus (see also Proposition
\ref{propendomor}) we get
$$
\pi(\delta(a))=\overline{\pi}([\delta(a),0,0,...])=\overline{\pi}(\Tdelta([a,0,0,...]))=
U\overline{\pi}([a,0,0,...])U^*=U\pi(a)U^*.
$$
$$
U^*U=\overline{\pi}(\chi_{\TDelta_{1}})=\overline{\pi}([\chi_{\Delta_1},0,0,...])=\pi(\chi_{\Delta_1}),
\,\,\,\,\,\,\,
UU^*=\overline{\pi}(\chi_{\TDelta_{-1}})=\overline{\pi}([\chi_{\Delta_-1},0,0,...])=\pi(\chi_{\Delta_{-1}}).
$$
Hence $(\pi,U,H)\in\mathrm{CovRep}\,(\A,\delta)$.\\
For $ [a_0,...,a_N,0,...]\in \B_0$ we have, cf.
Proposition \ref{propendomor},
$$
\overline{\pi}([a_0,...,a_N,0,...])=\overline{\pi}(a_0
+\Tdelta_*(a_1)+...+\Tdelta_*^N(a_N))=\pi(a_0)+U^*\pi(a_1)U+...+U^{*N}\pi(a_N)U^N.
$$
Hence the second part of the theorem  follows.
\end{Proof}
\begin{corl}
There is a natural bijection between
$\mathrm{CovFaithRep}\,(\A,\delta)$ and
$\mathrm{CovFaithRep}\,(\B,\Tdelta)$.
\end{corl}
The endomorphism $\Tdelta$ maps isomorphically  $C_{\TDelta_{-1}}(\X)$ onto $C_{\TDelta_{1}}(\X)$,
 whence we have a $*$-isomorphism between two closed two-sided
ideals  in $\B$.  In \cite{exel} R. Exel calls such  mappings partial automorphisms (of $\B$),
see \cite[Definition 3.1]{exel}. He also considers there
covariant representations of  partial automorphisms and his definition of those objects agrees
with Definition \ref{cov-rep-def} in case $\al$ is a partial homeomorphism. Moreover, R. Exel proves
in \cite[Theorem 5.2]{exel} the existance of covariant faithful representation of a partial automorphism
which automatically gives us
\begin{corl}\label{gupiecorl}
The set $\mathrm{CovFaithRep}\,(\A,\delta)$ is not empty.
\end{corl}
\begin{Proof}
The set $\mathrm{CovFaithRep}\,(\B,\Tdelta)$ is not empty by \cite[Theorem 5.2]{exel}.
\end{Proof}
The former of the preceding corollaries is not true for not faithful
representations (see Example \ref{simplexample2}). In general the
set $\mathrm{CovRep}\,(\B,\Tdelta)$ is larger than
$\mathrm{CovFaithRep}\,(\A,\delta)$, and there appears a problem
with prolongation of $\pi$ from $\A$ to $\B$ when
$\pi$ is not faithful. Fortunately in view of Proposition \ref{directlimit} we have the
following statement true.
\begin{thm}\label{coefalg3}
If $\delta:\A\rightarrow \A$ is a monomorphism then for any
$(\pi,U,H)\in\mathrm{CovRep}\,(\A,\delta)$ there exist
$(\overline{\pi},U,H)\in\mathrm{CovRep}\,(\B,\Tdelta)$
such that $\overline{\pi}([a_0,...,a_N,0,...])=
\pi(a_0)+U^*\pi(a_1)U+...+U^{*N}\pi(a_N)U^N$.
\end{thm}
\begin{Proof}
Let $(\pi,U,H)\in\mathrm{CovRep}\,(\A,\delta)$. Let us notice that
as $\delta$ is injective $\Delta_{-1}=X$ and hence $U$ is an
isometry (see Remark \ref{covreprem}). Now, consider the
$C^*$-algebra $C^*\big(\bigcup_{n=0}^\infty U^{*n} \pi(\A) U^n\big)$  and define the family of
mappings $\pi_n:\A_n\rightarrow C^*\big(\bigcup_{n=0}^\infty U^{*n} \pi(\A) U^n\big)$, $n\in
\N$, by the formula
$$
\pi_n(a)=U^{*n}\pi(a)U^n,\qquad a\in \A_n.
$$
Then the following diagram
$$
\begin{array}{ccc}
 \A_n & \,\,\,\stackrel{\delta_n}{\longrightarrow} &\A_{n+1}  \\

 & \,\,\, \searrow  \pi_{n}  &\,\,\,\,\,\,\,\,\, \downarrow \pi_{n+1}\\

 &       &      C^*\big(\bigcup_{n=0}^\infty U^{*n} \pi(\A) U^n\big)
\end{array}
$$
commutates. Hence, according to  Proposition
\ref{directlimit} there exists a unique $C^*$-morphism
$\overline{\pi}$ such that the diagram
$$
\begin{array}{ccc}
 \A_n & \,\,\,\stackrel{\Phi_{*n}\,\,}{\longrightarrow} & \B \\

 & \,\,\, \searrow  \pi_{n}  &\,\,\,\,\,\,\,\,\, \downarrow \overline{\pi}\\

 &       &      C^*\big(\bigcup_{n=0}^\infty U^{*n} \pi(\A) U^n\big)
\end{array}
$$
commutes. The hypotheses now follows.
\end{Proof}
\begin{corl}
If $\delta$ is a monomorphism then the mapping
$\overline{\pi}\rightarrow \overline{\pi}|_\A$ establishes a
bijection from $\mathrm{CovRep}\,(\B,\Tdelta)$
 onto $\mathrm{CovRep}\,(\A,\delta)$.
\end{corl}

In case $\delta$ is not a monomorphism two different covariant representations of
 $(\B,\Tdelta)$ may induce the same covariant representation of $(\A,\delta)$.
\begin{exm}\label{simplexample2}
   Not to look very far let us take the system $(X,\al)$ from Example \ref{simplexample}.
 The corresponding $C^*$-dynamical system is $(\A,\delta)$ where $\A=C(\{x_0,x_1,x_2,y_2\})\cong\mathbb{C}^4$ and
$\delta(a)=(0,a_{x_0},a_{x_1},a_{x_1})$  for $a=(a_{x_0},a_{x_1},a_{x_2},a_{y_2})\in \A$.
The coefficient algebra is
$$\B=C(\{\x_0,\x_1,\x_2,\y_0,\y_1,\y_2\})\cong\mathbb{C}^6$$
and  we check that, for $a\in \A$,
$$
[a,0,...]=(a_{x_0},a_{x_1},a_{x_2},a_{x_0},a_{x_1},a_{y_2}),\qquad
[0,a\delta(1),0,...]=(0,a_{x_0},a_{x_1},0,a_{x_0},a_{x_1}),$$
$$
[0,0,a\delta^2(1),0,...]=(0,0,a_{x_0},0,0,a_{x_0})\quad\mathrm{and}\quad
[0,0,...,0,a\delta^N(1),0...]\equiv 0,\,\,\mathrm{for}\,\, N>2.
$$
Moreover, we have
$\Tdelta(a_{\x_0},a_{\x_1},a_{\x_2},a_{\y_0},a_{\y_1},a_{\y_2})=(0,a_{\x_0},a_{\x_1},0,a_{\y_0},a_{\y_1})$.
It is now straightforward that
$(\overline{\pi}_1,U,\mathbb{C}^2)$ and
$(\overline{\pi}_2,U,\mathbb{C}^2)$ where $U=\left(\begin{array}{c
c}
 0   &     0  \\
 1   &     0  \\
\end{array}\right)$,
$$
\overline{\pi}_1(a_{\x_0},a_{\x_1},a_{\x_2},a_{\y_0},a_{\y_1},a_{\y_2})=\left(\begin{array}{c
c}
 a_{\x_0}   &      0      \\
 0          &    a_{\x_1}  \\
\end{array}\right)
\quad \mathrm{and} \quad
\overline{\pi}_2(a_{\x_0},a_{\x_1},a_{\x_2},a_{\y_0},a_{\y_1},a_{\y_2})=\left(\begin{array}{c
c}
 a_{\y_0}   &      0      \\
 0          &    a_{\y_1}  \\
\end{array}\right)
$$
are covariant representations of $(\B,\Tdelta)$
which induce the same covariant representation
$(\pi,U,\mathbb{C}^2)$ of $(\A,\delta)$ where
$\pi(a_{x_0},a_{x_1},a_{x_2},a_{y_2})=\left(\begin{array}{c c}
 a_{x_0}   &      0      \\
 0          &    a_{x_1}  \\
\end{array}\right)$.
\end{exm}

\section{Covariance
algebra}

In this section we introduce the title object of the paper. We
recall the definition of the partial crossed product, cf.
\cite{exel}, \cite{mcclanachan}, and then we define the covariance
algebra of $(\A,\delta)$  to be the  partial crossed product
associated with  $(\B,\Tdelta)$. We give a number of examples of
such algebras, and finally we show (justify the definition) that
the covariant algebra is the universal   object in the category of
covariant faithful representations  of $(\A,\delta)$ and in the
case $\delta$ is injective also in the category of covariant  (not
necessarily faithful) representations of $(\A,\delta)$.
\subsection{The
algebra $C^*(X, \al)$}

Let us  recall that  a \emph{partial automorphism} of a
$C^*$-algebra $ \mathcal{C}$ is a mapping $\theta:I\rightarrow J$
where $I$ and $J$ are closed two-sided ideals in $\mathcal{C}$ and
$\theta$ is a $*$-isomorphism, cf. \cite{exel}. If  a partial
automorphism $\theta$ is given then for each $n\in \Z$  we let
$D_n$ denote the domain of $\theta^{-n}$ with the convention that
$D_0=\mathcal{C}$ and $\theta^0$ is the identity automorphism of
$\mathcal{C}$. Letting then
$$\label{algebraL}
L=\{a\in l^1(\mathbb{Z},\mathcal{C}): a(n)\in D_n\}
$$
and defining  the convolution multiplication, involution and norm
 as follows
$$
(a*b)(n)=\sum_{k=-\infty}^{\infty}\theta^{k}\Big(\theta^{-k}\big(a(k)\big)b(n-k)\Big)
$$
$$
(a^*)(n)=\theta^n(a(-n)^*)
$$
$$
\|a\|=\sum_{n=-\infty}^{\infty}\|a(n)\|
$$
we equip $L$ with Banach $^*$-algebra structure. The universal
enveloping $C^*$-algebra of $L$ is called the \emph{partial
crossed product} (or the \emph{covariance algebra}) for the
partial automorphism $\theta$ and is denoted by
$\mathcal{C}\rtimes_{\theta} \mathbb{Z}$, see  \cite{exel},
\cite{mcclanachan}.
\\

It is clear  that the partial homeomorphism $\tal$ of $\X$ defines
the partial automorphism $\Tdelta:C_{\TDelta_{-1}}(\X)\rightarrow
C_{\TDelta_{1}}(\X)$ of the coefficient $C^*$-algebra $\B=C(\X)$
(we shall not distinguish  between the endomorphism $\Tdelta$
given by  \eqref{tildedelta} and its restriction to
$C_{\TDelta_{-1}}(\X)\subset \B$). The  definition  to follow
anticipates  Theorem \ref{UoCa}.
\begin{defn}\label{covariance_algebra}
The \emph{covariance algebra} $C^*(X,\al)$   of a partial dynamical system $(X,\al)$ is the partial crossed
 product for the partial automorphism $\Tdelta$ of the coefficient $C^*$-algebra $\B$.\\
That is $C^*(X,\al)=\B\rtimes_{\Tdelta}\Z$ and for $C^*(X,\al)$ we
shall  also write  $C^*(\A,\delta)$.
\end{defn}
\begin{rem} In the case
$\al$ is injective, equivalently $\delta$ is a partial
automorphism, the systems $(\A,\delta)$ and $(\B,\Tdelta)$
 are equal  and the covariance algebra of $(\A,\delta)$ is simply the partial crossed
 product. In particular, if $\al$ is a full homeomorphism, equivalently $\delta$ is an automorphism,
 then $C^*(\A,\delta)$ is the classic  crossed product.
 As we shall see, in the case $\al$ is surjective, equivalently $\delta$ is a monomorphism, $C^*(\A,\delta)$
 is the crossed product by a monomorphism considered for instance
 in \cite{Murphy}, \cite{Paschke}, \cite{exel2}, \cite{Adji_Laca_Nilsen_Raeburn}, cf. Remark \ref{murphyrem}.
\end{rem}
 Let $\sum_{k=-N}^{N}a_k u^k$ stands for  the element
 $a$ in $L=\{a\in l^1(\mathbb{Z},\B): a(n)\in C_{\TDelta_{n}}(\X)\}$ such that $a(k)=a_k$ for $|k|\leqslant N$,
and $a(k)=0$ otherwise. In view of the defined operations on $L$
it is clear  that $u$ is a partial isometry\label{u}, $u^k$ is $u$
to power $k$ and $(u^k)^*=u^{-k}$, so this  notation should not
cause any confusion.
  Using the natural injection $\B\ni a \mapsto au^{0}\in L$ we  identify
  $\B=C(\X)$
with the subalgebra of $C^*(\A,\delta)$, see \cite[Corollary 3.10]{exel}. Recalling the
identification  from Definition \ref{idenitify-defn} we have
$$
\A\subset \B\subset
C^*(\A,\delta)=\B\rtimes_{\theta}\Z.
$$
\begin{exm}\label{Toeplitz algebra}
The covariance algebra of the dynamical system considered in
Examples \ref{exm x_1-x_0} and \ref{exm x_1-x_0-coefalg} is
the Toeplitz algebra. Indeed, the coefficient algebra
$\B$ consists of convergent  and
$\Tdelta$ is a forward shift so the partial crossed product
$\B\rtimes_{\Tdelta}\Z$, cf. \cite{exel}, is the
Toeplitz algebra.
\end{exm}
\begin{exm}
Let us go back again to Example \ref{simplexample} (see also
Example \ref{simplexample2}). It is immediate to see that
$C^*(X,\al)=C^*(X',\al')$ and invoking \cite{exel}, or
\cite[Example 2.5]{mcclanachan} we can identify this algebra with
$M_3\oplus M_3$ where $M_3$ is the algebra of complex matrices
$3\times 3$. If we set $\A=C(X)$ and $\A'=C(X')$ then due to the
above remark we note that $\A$ and $\A'$ consist of the matrices
of the form
$$
\left(
\begin{array}{c c c}
 a_{x_0} &    0  &   0         \\
     0   & a_{x_1} & 0         \\
     0 &   0     &  a_{x_2}
 \end{array}\right)
  \oplus
\left(
\begin{array}{c c c}
 a_{x_0} &    0  &   0         \\
     0   & a_{x_1} & 0         \\
     0 &   0    &  a_{y_2}
 \end{array}\right),\qquad\mathrm{and} \qquad
\left(
\begin{array}{c c c}
 a_{x_0} &    0  &   0         \\
     0   & a_{x_1} & 0         \\
     0 &   0     &  a_{x_2}
 \end{array}\right)
  \oplus
\left(
\begin{array}{c c c}
 a_{x_0} &    0  &   0         \\
     0   & a_{y_1} & 0         \\
     0 &   0     &  a_{y_2}
 \end{array}\right)
  $$
respectively. The dynamics on $\A$ and $\A'$ are implemented by
the partial isometry
 $
U:=\left(
\begin{array}{c c c}
 0 &    0  &   0         \\
     1   & 0 & 0         \\
     0 &   1     &  0
 \end{array}\right)
  \oplus
\left(
\begin{array}{c c c}
 0 &    0  &   0         \\
    1   & 0 & 0         \\
     0 &  1    &  0
\end{array}\right)$.
The coefficient algebras $\B_0=\{\bigcup_{n=0}^\infty
U^{*n}(\A)U^n\}$ and $\B_0'=\{\bigcup_{n=0}^\infty
U^{*n}(\A')U^n\}$ equal with the algebra of diagonal matrices $
\left(
\begin{array}{c c c}
 a_{\x_0} &    0  &   0         \\
     0   & a_{\x_1} & 0         \\
     0 &   0     &  a_{\x_2}
 \end{array}\right)
  \oplus
\left(
\begin{array}{c c c}
 a_{\y_0} &    0  &   0         \\
     0   & a_{\y_1} & 0         \\
     0 &   0    &  a_{\y_2}
 \end{array}\right),
$ and it is straightforward  that $M_3\oplus
M_3=C^*(\A,U)=C^*(\A',U)$ and $\A\subset\A'$. In fact $\A$ is the
smallest $C^*$-subalgebra of $M_3\oplus M_3$ such that $U^*U\in
\A$,  $U(\cdot)U^*$ is an endomorphism of $\A$ and $C^*(\A,U)=M_3\oplus
M_3$.
 \end{exm}
\begin{exm}\label{finitedimension}
It is known, cf. \cite{exel},\cite{mcclanachan}, that an
 arbitrary finite dimensional $C^*$-algebra can be expressed as a covariance algebra
 (partial crossed product) of a certain dynamical system.
The foregoing example inspires us to present  the smallest such
system in the sense that space $X$ has the least number of points.
Let $\A=M_{n_1}\oplus ... \oplus M_{n_k}$ where $1\leqslant
n_1\leqslant ...\leqslant n_k$ be a finite dimensional
$C^*$-algebra and let us first assume that is there is no factor
$M_1$ in the decomposition of $\A$, that is $n_1>1$. We set
$X=\{x_1,x_2,...,x_{n_k-1},y_{n_1},y_{n_2},...,y_{n_k}\}$, so
$|X|=n_k+k-1$,
and $\al(x_m)=x_{m-1}$, for $m=2,...,n_k-1$; $\al(y_{n_m})=x_{n_m-1}$, for $m=1,...,k$,\\
\setlength{\unitlength}{.8mm}
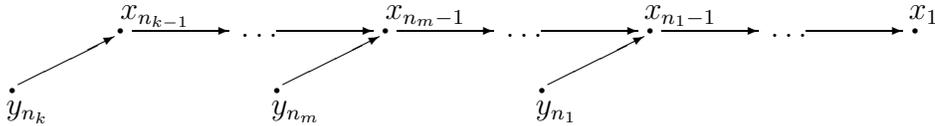
\begin{picture}(200,21)(-0,0)
\put(180,16){\circle*{1}} \put(179,18){$x_1$}
\put(156,15){$\dots$} \put(162,16){\vector(1,0){16}}
\put(138,16){\vector(1,0){16}}

\put(136,16){\circle*{1}} \put(135,18){$x_{n_1-1}$}
\put(119,7){\vector(2,1){16}} \put(118,6){\circle*{1}}
\put(117,2){$y_{n_1}$} \put(118,16){\vector(1,0){16}}
\put(112,15){$\dots$} \put(94,16){\vector(1,0){16}}

\put(92,16){\circle*{1}} \put(92,18){$x_{n_m-1}$}
\put(75,7){\vector(2,1){16}} \put(74,6){\circle*{1}}
\put(73,2){$y_{n_m}$} \put(74,16){\vector(1,0){16}}
\put(68,15){$\dots$} \put(50,16){\vector(1,0){16}}

\put(48,16){\circle*{1}} \put(48,18){$x_{n_{k-1}}$}
\put(31,7){\vector(2,1){16}} \put(30,6){\circle*{1}}
\put(29,2){$y_{n_k}$}
\end{picture}\\
 It is clear that
$C^*(X,\al)=\A$,  see \cite[Example 2.5]{mcclanachan}. In order to
include also algebras containing a number, say $l$, of  one-dimensional factors one
should simply add  $l$ points  to the above diagram.
\end{exm}
\begin{exm}
Let $a$ be the bilateral weighted shift on a separable Hilbert
space $H$ and let $a$ have the closed range. We have the polar
decomposition $a=U|a|$, where $|a|$ is a diagonal operator and $U$
is the bilateral shift.
 If we denote by $\A$ the commutative $C^*$-algebra $C^*(1,\{U^n|a|U^{*n}\}_{n\in \N})$,  then
 $\delta(\cdot)=U(\cdot)U^*$ is a unital injective  endomorphism of $\A$ and hence the dynamical system
 $(X,\al)$ corresponding to $(\A,\delta)$ is  such that $\al:X\rightarrow X$ is onto.
 It is immediate that the coefficient algebra $\B=C(\X)$ has the form $C^*(1,\{U^n|a|U^{*n}\}_{n\in \Z})$.
  Thus due to \cite[Theorem 2.2.1]{Odonovan},
 $C^*(a)=\B\rtimes_{\Tdelta} \Z$ and  so
$$
C^*(a)=C^*(X,\al).
$$
Following \cite{Odonovan} we present now  the canonical form of
$(X,\al)$. Let $Y$ denote the spectrum of $|a|$ and  let
$T:X\rightarrow \prod_{n=0}^\infty Y$ be defined by
$T(x)=(x(|a|),x(\delta(|a|)),...,x(\delta^n(|a|)),... )$. Then
similarly to \cite{Odonovan} we infer that $T$ is a homeomorphism
of $X$ onto $T(X)$, where  $T(X)$ is given the topology induced by
the product topology on $\prod_{n=0}^\infty Y$, and under $T$,
$\al$ becomes a  shift on the  product space  $T(X)$.
\end{exm}

\subsection{Universality of $C^*(X,\al)$ }
 Now we
are in  position  to prove the main result of this section
which justifies the anticipating Definition
\ref{covariance_algebra}. We shall
base on the results from the previous section and some known facts
concerning the partial crossed product \cite{exel}.
 We adopt  to the commonly used notation
$U^{-n}=U^{*n}$ where  $U$ is a partial isometry and $n\in\N$.

\begin{thm}\label{UoCa} Let
$(\pi,U,H)\in\mathrm{CovFaithRep}\,(\A, \delta)$ or
$(\pi,U,H)\in\mathrm{CovRep}\,(\A, \delta)$ in the case $\delta$ is a
monomorphism. Then the formula
\begin{equation}\label{kopytko}
(\overline{\pi}\times U)(\sum_{n=-N}^{N}a(n) u^n)=\sum_{n=-N}^N
\Big(\sum_{k=0}^\infty U^{*k}\pi(a_k^{(n)})U^k\Big)U^n
\end{equation}
where $a(n)=[a_0^{(n)},a_1^{(n)},...,a_k^{(n)},..]
\in\B_0$,  establishes   an epimorphism  of the
covariance algebra $C^*(\A,\delta)$ onto the $C^*$-algebra
$C^*(\pi(\A),U)$ generated by $\pi(\A)$ and $U$.
\end{thm}
\begin{Proof}
In both cases, $(\pi,U,H)\in\mathrm{CovFaithRep}\,(\A, \delta)$ or
$(\pi,U,H)\in\mathrm{CovRep}\,(\A, \delta)$ and $\delta$ is
injective, $(\pi,U,H)$ extends to the   covariant representation
$(\overline{\pi},U,H)$ of $(\B,\Tdelta)$, see
Theorems \ref{coefalgebra2}, \ref{coefalg3}. Since $C^*(\pi(\A),U)=\B\rtimes_{\Tdelta}\Z$ we
obtain,  by  \cite[Proposition 5.5]{exel}, that
$$
(\overline{\pi}\times
U)(\sum_{n=-N}^{N}a(n) u^n)=\sum_{n=-N}^N\overline{\pi}(a(n)) U^n.
$$
establishes  the representation of
$C^*(\A,\delta)$, and by Theorems \ref{coefalgebra2} and \ref{coefalg3}, $(\overline{\pi}\times
U)$ is in fact given by \eqref{kopytko}.
\end{Proof}

Due to Definition \ref{cov-rep-def}, Corollary \ref{gupiecorl}, and the preceding
Theorem \ref{UoCa} we can  alternatively define the covariance algebra
to be the  universal unital $C^*$-algebra
generated by a copy of $\A$ and a partial isometry $u$ subject to
relations
$$
  u^*u\in \A, \qquad \quad \delta(a) = uau^* ,\qquad a\in \A,
$$
see also Proposition \ref{prop(O)sition}. In particular the above
relations imply that $uu^*=P_{\Delta_1}$ and
$u^*u=P_{\Delta_{-1}}$. Thus  $\delta$ is injective iff $u$ is
isometry, and $\delta$ is an automorphism iff $u$ is unitary.
\begin{thm}\label{uocainvers} Let $\sigma$ be a representation of $C^*(\A, \delta)$
 on a Hilbert space $H$.
  Let  $\pi$ denote the restriction of $\sigma$ onto $\A$ and let  $U=\sigma(u)$. Then
   $(\pi,U,H)\in \mathrm{CovRep}\,(\A,\delta)$.
\end{thm}
\begin{Proof}
By \cite[Theorem 5.6]{exel} we have $(\overline{\pi},U,H)\in
\mathrm{CovRep}\,(\B,\Tdelta)$ where
$\overline{\pi}$ is $\sigma$ restricted to $\B$.
Hence by Theorem \ref{dobrethm} we get $(\pi,U,H)\in
\mathrm{CovRep}\,(\A,\delta)$.
\end{Proof}

\begin{corl}\label{murpytheorem} If $\delta$ is a monomorphism then the correspondance
$(\pi,U,H) \longleftrightarrow (\overline{\pi}\times U)$, cf. Theorem \ref{UoCa}, is a bijection
between $\mathrm{CovRep}\,(\A, \delta)$ and the set
of all representations of $C^*(\A,\delta)$.
\end{corl}
\begin{Proof}
In virtue of  Theorem \ref{UoCa} the mapping $(\pi,U,H)
\longleftrightarrow (\overline{\pi}\times U)$ is a well defined
injection and by Theorem \ref{uocainvers} it is also a surjection.
\end{Proof}
\begin{rem}\label{murphyrem} Corollary \ref{murpytheorem}  can be considered as a special case
of Theorem 2.3  from the paper \cite{Murphy} where
(twisted) crossed products by injective endomorphisms were
investigated. However, our approach is slightly different.
Oversimplifying; Murphy  defines the algebra $C^*(\A,\delta)$ as
$pZp$ where $Z$ is the full crossed product of direct limit
$B=\underrightarrow{\lim\,\,\,}\A$ and $p$ is a certain projection
from $B$, whereas we include the projection $p$ in the direct
limit $\B=\underrightarrow{\lim\,\,\,}\A_n \subset
B$, see Proposition \ref{directlimit},  and hence take the partial
crossed product.
\end{rem}
\section{Invariant subsets and the topological freedom of  partial
mappings}

The present section is  devoted to the generalization of two
important
 notions of  the theory of   crossed products.
  We start with  a definition of $\al$-invariant sets.  With  help of this notion we
  will   describe (in the next section)
   the ideal structure
   of
 the  covariance algebra.
   Next we introduce  a definition of topological freedom - a property which is a powerful instrument when
    used to construct faithful representations of covariance algebra.
    Thanks to that in Section 6 we prove a  version of The Isomorphism Theorem.
\subsection{Definition of the  invariance under $\al$ and $\tal$ and their interrelationship}
   The  Definition \ref{alphainvar} to follow might look  strange at first however  the author's
   impression is
   that among the others this one is the most natural  generalization of that from
   \cite[Definition 2.7]{exel_laca_quigg} for
   the case considered here. One may treat  $\al$-invariance as  the invariance under
   the partial action of $\mathbb{N}$ on $X$.
\begin{defn}\label{alphainvar}
Let $\al$ be a  partial mapping  of $X$. A subset $V$ of $X$ is
said to be \emph{invariant} under the partial mapping $\al$, or
shorter \emph{$\al$-invariant}, if
\begin{equation}\label{defninvareq} \al^n(V\cap \Delta_n)= V\cap
\Delta_{-n},\qquad n=0,1,2...\,. \end{equation}
\end{defn}
 When $\al$ is  injective then we have another   mapping
 $\al^{-1}:\Delta_{-1}\rightarrow X$ and life is a bit easier.
\begin{prop}\label{alphainvariance}
Let $\al$ be an injective partial mapping and let $V\subset X$.
Then  $V$ is $\al$-invariant if and only if one of the following
conditions holds
\begin{description}
\item[i)]$V$ is $\al^{-1}$-invariant
 \item[ii)] for each $n=0, \pm
 1, \pm 2,...,
 $ we have $\al^n(V\cap\Delta_{n})\subset V$
\item[iii)] $\al(V\cap\Delta_{1})\subset V\qquad \mathrm{and}
\qquad \al^{-1}(V\cap\Delta_{-1})\subset V$,
 \item[iv)]  $\al(V\cap \Delta_1)= V\cap \Delta_{-1}$.
\end{description}
\end{prop}
\begin{Proof}
The equivalence of invariance of $V$ under  $\al$ and $\al^{-1}$
is straightforward, and so are   implications i)$\Rightarrow$
ii)$\Rightarrow$ iii).
 \\
 To prove iii)
$\Rightarrow$ iv) let us observe  that since
$\al(V\cap\Delta_1)\subset \Delta_{-1}$ and
$\al^{-1}(V\cap\Delta_{-1})\subset \Delta_{1}$ we get
$\al(V\cap\Delta_{1})\subset V\cap \Delta_{-1}$ and
$\al^{-1}(V\cap\Delta_{-1})\subset V\cap \Delta_{1}$. The latter relation implies that  $V\cap
\Delta_{-1}=\al(\al^{-1}(V\cap
\Delta_{-1}))\subset\al(V\cap\Delta_{1})$ and so $\al(V\cap
\Delta_1)= V\cap \Delta_{-1}$.
\\
The only thing   left to  be shown is that iv) implies
$\al$-invariance of $V$. We prove this by  induction.
Let us assume that $\al^k(V\cap \Delta_k)= V\cap \Delta_{-k},$ for
$k=0,1,...,n-1\,$. By injectivity it is equivalent to $
\al^{-k}(V\cap \Delta_{-k})= V\cap \Delta_{k},$ for
$k=0,1,...,n-1\,$. As $\Delta_n\subset \Delta_{n-1}$ and
$\Delta_{-n}\subset \Delta_{-(n-1)}$ we have $
 \al^n(V\cap \Delta_n)\subset V\cap \Delta_{-n}$ and $\al^{-n}(V\cap \Delta_{-n})\subset V\cap
\Delta_{n}. $ Applying $\al^n$ to the latter relation we get
$V\cap \Delta_{-n} \subset \al^n(V\cap \Delta_n)$ and hence $
\al^n(V\cap \Delta_n)= V\cap \Delta_{-n}$.
\end{Proof}
 Item ii) tells us  that Definition \ref{alphainvar} extends
   \cite[Definition 2.7]{exel_laca_quigg}  in the case of a single partial mapping.
  Let us note that if $\Delta_1\neq X$ and $\al$ is not injective, then none of items ii)-iv) is
  equivalent to (\ref{defninvareq}) and therefore none of them could be used as a definition of $\al$-invariance.
\begin{exm}\label{exm alfainvar}
Indeed, let  $X=\{x_0,x_1,x_2,y_2,y_3\}$  and
let $\al$ be defined by the
relations $\al(y_3)=y_2$, $\al(y_2)=\al(x_2)=x_1$ and  $\al(x_1)=x_0$:\\
\setlength{\unitlength}{.8mm}
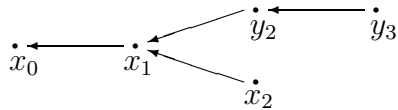
\begin{picture}(80,20)(-65,0)
\put(15,10){\circle*{1}} \put(14,6){$x_0$}
\put(35,10){\circle*{1}} \put(33,6){$x_1$}
\put(55,16){\circle*{1}} \put(54,12){$y_2$}
\put(75,16){\circle*{1}} \put(74,12){$y_3$}
\put(55,4){\circle*{1}} \put(53,0.5){$x_2$}
\put(53,16){\vector(-3,-1){16}} \put(33,10){\vector(-1,0){16}}
\put(53,4){\vector(-3,1){16}} \put(73,16){\vector(-1,0){16}}
\end{picture}\\
 Then the set $V_1=\{x_0,x_1,x_2\}$ fulfills
item iv)  but it is not invariant under $\al$ in the sense of
Definition \ref{alphainvar}. Whereas the set
$V_2=\{x_0,x_1,y_2,y_3\}$ is $\al$-invariant but it does not
fulfill items ii) and iii).
\end{exm}

 We shall show that in general there are less $\al$-invariant sets  in $X$  than $\tal$-invariant sets in
 $\X$, cf. Theorem \ref{latticeisomorph} and a remark below,
but fortunately there is  a natural  one-to-one correspondence between closed
  sets  invariant under $\al$ and closed sets invariant under $
  \tal$, see Theorem \ref{latticeiso}.   We start with  a lemma.
\begin{lem}\label{invarianceprop}
Let $\al$ be a partial mapping of $X$ and let $U\subset X$ be
invariant under $\al$. Then we have
$$\al^k(\Delta_{n+k}\cap U)=U\cap \Delta_n\cap \Delta_{-k},\qquad
n,k=0,1,2,...\,.$$
\end{lem}
\begin{Proof}
As $\al^k(U\cap\Delta_k)=U\cap \Delta_{-k}$ and
$U\cap\Delta_{n+k}\subset U\cap \Delta_{k}$ we have
$\al^k(\Delta_{n+k}\cap U)\subset U\cap \Delta_n\cap \Delta_{-k}$,
for $k,n\in \N$. On the other hand,  for every $x\in U\cap
\Delta_n\cap \Delta_{-k}$ there exists $y\in U\cap\Delta_k$ such
that $\al^k(y)=x$. Since $x\in \Delta_{n}$ we we have $y\in
U\cap\Delta_{n+k}$ and hence $U\cap \Delta_n\cap \Delta_{-k}
\subset \al^k(\Delta_{n+k}\cap U)$.
 \end{Proof}
\begin{thm}\label{latticeisomorph}
Let $(\al, X)$ be a partial dynamical system and let $(\tal,
\X)$ be its reversible extension. Let $\Phi:\X\rightarrow X$ be
the projection defined by (\ref{Phimap}). Then  the map
\begin{equation}\label{latticemap}
 \X\supset \U\longrightarrow
U=\Phi(\U)\subset X
\end{equation}
 is a surjection from the family of
$\tal$-invariant subsets of $\X$ onto  the family of
$\al$-invariant subsets of $X$. Furthermore if $\{U_n\}_{n\in\N}$
are the sections  of $\U$ (see Definition \ref{U_n})), then $\U$ is
$\tal$-invariant
 if and only if $U_0$ is $\al$-invariant and
 $$ U_n=U_0 \cap \Delta_n,\qquad n=0,1,2,...\, .$$
\end{thm}
\begin{Proof}
Let $\U$ be $\tal$-invariant and let $U=\Phi(\U)$. Then by
(\ref{tildealpha}) and $\tal$-invariance of $\U$, for each $n\in
\N$, we have
$$\al^n(U\cap \Delta_n)=\Phi(\tal^n(\U\cap\TDelta_n))=
 \Phi(\U\cap\TDelta_{-n})=U\cap
 \Delta_{-n}$$
and hence $U$ is invariant under $\al$ and  the mapping
(\ref{latticemap}) is well defined. Moreover, if $U_n$,
$n\in\mathbb{N}$, are the sections of $\U$, then by invariance of
$\U$ under $\tal^{-1}$ (see Proposition \ref{alphainvariance}) we
get
$$
U_n=\Phi(\tal^{-n}(\U\cap \TDelta_{-n})) =
\Phi(\U\cap\TDelta_n)=U\cap\Delta_n
$$
where $U=U_0$ is $\al$-invariant. \\
Now, we show that the mapping (\ref{latticemap}) is surjective.
Let $U$ be any nonempty set invariant under $\al$ and let us
consider $\U$ of the form
$$
\U:=\Big(U\times (U\cap \Delta_1\cup\{0\})\times ...\times
(U\cap\Delta_n \cup\{0\})\times ...\Big)\cap\X.
$$
What we need to prove is that  $\Phi(\U)=U$ (note that in
general  $\U$ may  occur   to be empty). It is clear that
$\Phi(\U)\subset U$. In order to prove that  $U \subset\Phi(\U)$
we  fix an arbitrary point $x_0\in U$ and  suppose that there does not exist  $\x$ in $\U$
such that $\Phi(\x)=x_0$. We will  construct an infinite sequence $(x_0,x_1,x_2,...)$ in $\U$
and thereby obtain a contradiction.
\\
Indeed, we must have   $x_0\in U\cap \Delta_{-1}$, for otherwise we can take $\x=(x_0,0,0,...)\in \U$.
 Hence  by Lemma \ref{invarianceprop} there exists
$x_1\in U\cap\Delta_1$ such that $\al(x_1)=x_0$.
Suppose now  we have chosen  $n-1$ points $x_1,...,x_n$ such that $x_k\in U\cap \Delta_{k}$
and $\al(x_{k})=x_{k-1}$ for $k=1,...,n$, then
  $x_{n}$ must be in $U\cap\Delta_{n}\cap\Delta_{-1}$, for otherwise  we can take $\x=(x_0,x_1,...,x_{n},0,...)\in
\U$. Hence by Lemma \ref{invarianceprop} there exists
$x_{n+1}\in U\cap\Delta_{n+1}$ such that $\al(x_{n+1})=x_n$. This  ensure us that there is a
 sequence $\x=(x_0,x_1,x_2,...)$ such that $x_n\in U\cap
\Delta_n$ and $ \al(x_n)=x_{n-1} $, for all $n=1,2,...\,$. Thus $\x\in\U$ and we arrive at the contradiction.
\\
 By virtue of item v) in Proposition \ref{alphainvariance} in order to prove
the  $\tal$-invariance of $\U$ it suffices
 to show that
$$
\tal(\U\cap\TDelta_{1})\subset \U\qquad \mathrm{and} \qquad
\tal^{-1}(\U\cap\TDelta_{-1})\subset \U
$$
but this follows  immediately  from the form of $\U$, $\tal$,
$\tal^{-1}$ and from $\al$-invariance of $U$. Thus according to
the first part of the proof we conclude that, for each $n\in \N$,
the $n$-section  $U_n$ of $\U$ is equal to $U\cap \Delta_n$. The proof is complete.
\end{Proof}
\begin{corl}\label{latticeiso1}
If $\al$ is injective on the inverse image of
$\Delta_{-\infty}=\bigcap_{n\in\N}\Delta_{-n}$ (for  $x\in
\Delta_{-\infty}$ we have $|\al^{-1}(x)|=1$) then the mapping
(\ref{latticemap})  from the previous theorem is a bijection.
\end{corl}
\begin{Proof} It suffices to apply Proposition \ref{Implication}.
\end{Proof}
 Under the hypotheses of Theorem \ref{latticeisomorph}, surjection considered
there  might not be  a bijection. For instance in Example
\ref{exm x_1-x_0} we have three $\al$-invariant sets:
$X,\{x_0\},\emptyset$, and four
 $\tal$-invariant sets: $\overline{\N},\N,\infty,\emptyset$.
However the mapping \eqref{latticemap} is always bijective when
restricted to closed invariant sets.
\begin{thm}\label{latticeiso}
The mapping (\ref{latticemap})  is a bijection  from the family of
$\tal$-invariant closed sets onto the family of $\al$-invariant
closed sets.
\end{thm}
\begin{Proof}
 By Theorems
\ref{latticeisomorph} and \ref{closed set brackets}, for every
 $\tal$-invariant  closed subset $\U$ such that
$\Phi(\U)=U$ we have
\begin{equation}\label{Ubrackets}
\U=\Big(U\times (U\cap \Delta_1\cup\{0\})\times ...\times
(U\cap\Delta_n \cup\{0\})\times ...\Big)\cap\X,
\end{equation}
that means $\U$ is uniquely determined by $U$. Since $\X$ and $X$
are compact, and $\Phi:\X\rightarrow X$ is continuous the set $U$
is closed. Hence $\Phi$ maps injectively the family of closed
$\tal$-invariant sets into the family of closed
$\al$-invariant sets.\\
 On the other hand, if $U$ is closed and $\al$-invariant then
by definition of the product topology, the set $\U$ given by
(\ref{Ubrackets}) is also closed and according to  Theorem
\ref{latticeisomorph}, $\Phi(\U)=U$. Thus the proof is complete.
\end{Proof}
The next important notion which we shall need  to obtain   a simplicity
criteria for covariance algebra (see Corollary \ref{simplicity}) is the notion of minimality, cf.
\cite{exel_laca_quigg}.
\begin{defn}
A partial continuous mapping  $\al$ (or a partial dynamical
system $(X,\al)$) is said to be \emph{minimal} if there are no
$\al$-invariant closed subsets of $X$ other than $\emptyset$ and
$X$.
\end{defn}
\begin{prop}\label{latticeiso2}
A partial dynamical system $(X,\al)$ is minimal if and only if its
reversible extension $(\X,\tal)$ is minimal.
\end{prop}
\begin{Proof} An easy consequence of Theorem \ref{latticeiso}.
\end{Proof}
When $\al$ is injective, the binary operations $"\cup"$ and
$"\cap"$, or equivalently partial order relation $"\subset"$,
define the lattice structure on the family  of $\al$-invariant
sets, see \cite{exel_laca_quigg}.
 The situation changes when $\al$ is not
injective. Of course $"\subset"$ is still a partial order relation
which determines the lattice structure, but it may happen that the
intersection of two $\al$-invariant sets is no longer
$\al$-invariant.
\begin{exm} Let $(X,\al)$ and $(X',\al')$ be  dynamical systems from  Example \ref{simplexample}.
It is easy to verify that  there are four
 sets invariant under $\al$: $X$, $V_1=\{x_0,x_1,x_2\}$, $V_2=\{x_0,x_1,y_2\}$ and $\emptyset$;
 and four sets invariant under $\al'$: $X'$, $V_1'=\{x_0',x_1',x_2'\}$, $V_2'=\{x_0',y_1',y_2'\}$,  $\emptyset$.
 Hence neither  $V_1\cap V_2$ nor $V_1'\cap V_2'$ is   invariant. However there are four invariant subsets: $\X$,
$\V_1=\{\x_0,\x_1,\x_2\}$, $\V_2=\{\y_0,\y_1,\y_2\}$, $\emptyset$
on the reversible extension level
 ($(\X,\tal)= (\X',\tal')$)
and  $\V_1\cap \V_2=\emptyset$ is invariant of course.
\end{exm}
\begin{defn}
We denote by $clos_\al(X)$ the lattice of $\al$-invariant closed
subsets of $X$ where the lattice structure is defined by the
partial order relation $"\subset"$.
\end{defn}
According to Theorem \ref{latticeiso}, $\Phi$ determines the lattice isomorphism $ clos_{\tal}(\X) \cong clos_\al(X)$.
\subsection{Topological freedom}
Recall now that  a partial action of a group $G$ on  a topological  space
$X$ is said to be topologically free if the set of fixed points
$F_t$, for each partial homeomorphism $\al_t$ with $t\neq e$, has
empty interior \cite{exel_laca_quigg}. In view of that, the next
definition constitutes a generalization of topological freedom
notion to the class of systems where dynamics are implemented by
one, not necessarily injective, partial mapping.
\begin{defn}\label{topological freedom}
Let  $\al:\Delta\rightarrow X$ be a  continuous partial mapping of
Hausdorff's topological space $X$. For each $n>0$, we set
$F_n=\{x\in\Delta_n: \al^n(x)=x\}$. It is said that the action of
$\al$ (or briefly $\al$) is \emph{topologically free}, if every
open nonempty subset $ U \subset F_n$ possess 'an exit', that is
there exists a point $x\in U$ such that one of the equivalent
conditions hold
 \Item{i)} for some $k=1,2,...,n$ we
have $ |(\al^{-k}(x))|>1,$
 \Item{ii)} for some  $k=1,2,...,n$
$\al^{-1}(\al^k(x))\neq \{\al^{k-1}(x)\}$.
\end{defn}
We supply now some   characteristics of this  topological freedom
notion.
\begin{prop}
The following conditions are equivalent
\begin{description}
\item[i)] $\al$ is topologically free, \item[ii)] for each $n>0$
and every open nonempty subset $ U \subset F_n$ there exist points
$x\in U$, $y\in \Delta_1\setminus F_n$ and a number $k=1,2,...,n$,
such that
$$
\al(y)=\al^{k}(x).
$$
\item[iii)] for each $n>0$, the set
$$
\{x\in \Delta_{n-1}: \al^{k-n}(x)=\{\al^k(x)\} \,\,\mathrm{for\,\,
}\,\, k=0,1,...,n-1 \}
$$
has an empty interior.
\end{description}
\end{prop}
\begin{Proof}
$i)\Rightarrow ii)$. Let $U\subset F_n$ be an open nonempty set.
Let
 $x\in U$ and $k=1,...,n$ be such that  item ii) from Definition \ref{topological
freedom} holds. We take $y\in\al^{-1}(\al^k(x))$ such that $y\neq
\al^{k-1}(x)$. Then $\al(y)=\al^k(x)$ and since
$\al^n(y)=\al^{k-1}(x)$ we have
$y\notin F_n$.\\
$ii)\Rightarrow iii)$.  Suppose  that for some $n>0$ there exists a
nonempty open subset $U$ of $ \{x\in \Delta_{n-1}:
\al^{k-n}(x)=\{\al^k(x)\} \,\,\mathrm{for\,\, }\,\, k=0,1,...,n-1
\} $. It is clear that $U\subset F_n$ and hence for some
$k_0=1,2,...,n$, there exists $y\in \Delta_1\setminus F_n$ such
that $\al(y)= \al^{k_0}(x)$.  Taking
$k=k_0-1$ we obtain that $y\in\al^{k-n}(x)=\{\al^{k}(x)\}$ and thus we arrive at a contradiction since
 $\al^{k}(x)\in F_n$ and $y\notin F_n$.\\
$iii)\Rightarrow i)$. Suppose that $\al$ is not topologically
free. Then there exists an open nonempty set $U\subset F_n$ such
that for all $x\in U$ and $k=1,...,n$, we have $|\al^{-k}(x)|=1$.
It is not hard to see that $U\subset \{x\in \Delta_{n-1}:
\al^{k-n}(x)=\{\al^k(x)\} \,\,\mathrm{for\,\, }\,\, k=0,1,...,n-1
\}$ and thereby we arrive at the contradiction. \end{Proof}
 The role similar to the one
which  topological freedom plays in the theory of crossed
products is the role played  in the theory of $C^*$-algebras associated with graphs
by the condition that every circuit  in a graph has an exit
(see, for example, \cite{pask}, \cite{exel_laca}). The connection
between these two properties is not only of theoretical character, see for
instance  \cite[Proposition 12.2]{exel_laca}).
  Taking this into account the  following two simple examples might be of interest.
   Before that let us establish the indispensable notation.
   \\
Let $A=(A(i,j))_{i,j\in\{1,...,N\}}$ be the matrix with entries in
$\{0,1\}$. It can be regarded as an \emph{adjacency matrix} of a
directed graph $Gr(A)$: the  vertices of $Gr(A)$ are numbers
1,...,N and  edges are pairs $(x,y)$ of vertices such that
$A(x,y)=1$. By a \emph{path} in $Gr(A)$ we mean a  sequence
$(x_0,x_1,...x_n)$ of vertices such that $A(x_k,x_{k+1})=1$ for
all $k$. A \emph{circuit}, or a \emph{loop}, is a finite path
$(x_0,...,x_n)$ such that $A(x_n,x_0)=1$. Finally a circuit
$(x_0,x_1,...,x_n)$ is said to have an \emph{exit} if, for some
$k$, there exists $y\in\{1,...,N\}$ with $A(x_k,y)=1$ and $y\neq
x_{k+1\,(mod\,\, n)}$.
\begin{exm}\label{graph}
Let $(X,\al)$ be a partial dynamical system such that
$X=\{1,...,N\}$ is finite. If we define $A$ by the relation:
$A(x,y)=1$ iff $\al(y)=x$, then $\al$ is
topologically free if and only if every loop in $Gr(A)$ has an
exit. This is an easy consequence of the fact that $Gr(A)$ is the
graph of the partial mapping $\al$ with reversed edges.
\end{exm}
 We shall say that a circuit $(x_0,x_1,...,x_n)$  has an \emph{entry} if,
 for some $k$, there exists $y\in\{1,...,N\}$ with $A(y,x_k)=1$ and $y\neq x_{k-1}\,(mod\,\, n)$.

\begin{exm}
Let $(X_A,\sigma_A)$ be a dynamical system where $\sigma_A$ is a
one-sided Markov subshift associated with a  matrix $A$, see page
\pageref{one-sided Markov}.  One-sided subshift $\sigma_A$ acts
topologically free if and only if every circuit in $Gr(A)$ has an
exit or an entry. Indeed, if  there exist a loop $(y_0,...,y_n)$
in $Gr(A)$ which has no exit  and no entry then $U=\{(x_k)_{k\in
\N}\in X_A: x_0=y_0\}=\{(y_0,y_1,...,y_n,y_0,...)\}$ is an open singleton, $U\subset F_{n+1}$ and $U$  has
no 'exit' in the sense of Definition \ref{topological freedom}, thereby $\sigma_A$ is not topologically free.
\\
On the
other hand, if  every loop in $Gr(A)$ has an exit or an entry, and if $x=(x_1,...,x_n,x_1,...)$
is an element of an open subset $U
\subset F_n$ for some $n>0$, then  the loop $(x_1,...,x_n)$ must have an entry, as  it clearly has  no exit.
Hence $(x_1,...,x_n)$   we have $
|(\sigma_A^{-k}(x))|>1$, for some $k=1,2,...,n$, and thus $\sigma_A$ is topologically free.
\end{exm}
We end this section with the result which, in a sense, justifies
Definition \ref{topological freedom}, and is the main tool
used to prove the Isomorphism Theorem.
\begin{thm}\label{topologicalfreedom}
Let $F_n=\{x\in\Delta_n: \al^n(x)=x\}$ and $\F_n=\{\x\in\TDelta_n:
\tal^n(\x)=\x\}$, $n\in\mathbb{N}\setminus\{0\}$. We have
\begin{equation}
\label{F_nset} \F_n=\{(x_0,x_1,...)\in \X: x_k\in
F_n,\,\,k\in\N\},\qquad n=1,2,...,
\end{equation}
 and  $\al$ is
topologically free if and only if  $\tal$ is topologically free.
\end{thm}

\begin{Proof}
Throughout the proof we fix an $n>0$. It is clear that $F_n$ and
$\F_n$ are invariant under $\al$ and $\tal$ respectively (see Definition \ref{alphainvar}), and that
$ \Phi(\F_n)=F_n$. By virtue of  Theorem \ref{latticeisomorph}
we have
$$
\F_n=\Big(F_n\times (F_n\cap \Delta_1\cup\{0\})\times ...\times
(F_n\cap\Delta_k \cup\{0\})\times ...\Big)\cap\X.
$$
But, since $F_n\subset \bigcap_{k\in\mathbb{Z}}\Delta_{k}$  we
obtain  $\F_n=(F_n\times F_n\times ...\times
F_n\times ...)\cap\X$ and hence (\ref{F_nset}) holds.\\
Now suppose that $\alpha$ is topologically free and on the
contrary that there exists an open nonempty subset $\U\subset
\F_n$. Without loss of generality, we can assume that it has the
form
$$\U=\Big(U_0\times U_1  ...\times
U_m\times \Delta_{m+1}\cup\{0\}  \times \Delta_{m+2}\cup\{0\} \times ...\Big)\cap\X
$$
where $U_0,U_1,...,U_m$ are open subsets of $X$, and as $\U\subset  \F_n=(F_n\times F_n\times ...\times
F_n\times ...)\cap\X$,  they are in fact subsets of $F_n$, and it is readily checked that
$$
\U=\Big(F_n \times F_n   ...\times \bigcap_{k=0}^m
\al^{-k}(U_{m-k})\times \Delta_{m+1}\cup\{0\}  \times \Delta_{m+2}\cup\{0\} \times ...\Big)\cap\X.
$$
The set $U:= \bigcap_{k=0}^m
\al^{-k}(U_{m-k})$ is an open and nonempty subset of $F_n$. Hence, due to
the topological freedom of $\al$ there exists $y\notin F_n$ and
$k=1,...,n$, such that $\al(y)=\al^k(x)$ for some $x\in U$. Taking any element $\x=(x_0,x_1,...)\in \X$  such that
$x_m:=x$,  $x_{m+i}:=\al^{n-i}(x)$
for $i=1,...,n-k$, and $x_{m+n-k+1}=y$ we arrive at the
contradiction, because  $\x\in \U$ and $\x \notin \F_n$.
\\
Finally suppose that $\alpha$ is not topologically free. Then
there exists an open nonempty subset $U\subset F_n$ such that, for
all $x\in U$, $|\al^{-k}(x)|=1$ and so
$$
\U=\{(x,\al^{n-1}(x),\al^{n-2}(x),...,\al^{1}(x),x,\al^{n-1}(x),...
)\in \X: x\in U\}=(U\times (X\cup \{0\}) \times ...)\cap \X
$$
is an open nonempty subset of $\F_n$. Hence $\tal$ is not
topologically free and the proof is complete.
\end{Proof}

\section{Ideal structure of covariance algebra and 
\\the Isomorphism Theorem}
It is well-known that every closed ideal of $\A=C(X)$ is of the
form $C_U(X)$ where $U\subset X$ is open, and therefore we have an
order preserving bijection between open sets and ideals. The
Theorem 3.5 from \cite{exel_laca_quigg} can be regarded as a
generalization of this fact; it says that, under some assumptions,
there exists a lattice isomorphism between open invariant sets and
ideals of the partial crossed product. In this section we shall
prove the new useful variant of this theorem. The novelty is that in our
approach (cf. Theorem \ref{latticeiso}) it is  more natural to
investigate a correspondence between ideals of the covariance
algebra and closed invariant sets. \\
After that we shall prove  the main result of this paper, a
version of the Isomorphism Theorem  where the main achievement is
that we do not assume any kind of reversibility  of an action on a
spectrum of a $C^*$-dynamical system. \par
\subsection{Lattice isomorphism of closed $\al$-invariant sets onto ideals of $C^*(X,\al)$}
 Let us start with the
 proposition which is an attempt of describing the concept of invariance on the algebraic level, cf.
 \cite[Definition 2.7]{exel_laca_quigg}. For that purpose we will abuse notation concerning endomorphism
 $\delta$ and denote by $\delta^n$, $n\in \N$, morphisms
$
 \delta^n:C(\Delta_{-n})\rightarrow C(\Delta_{n})
$
 of composition with  $\al^n:\Delta_n\rightarrow \Delta_{-n}$. We believe that this notation does not cause confusion,
  although we stress that set $\Delta_{-n}$ does not have to be open and hence we can not identify $C(\Delta_{-n})$ with
   a subset of $C(X)$. For instance, it may happen that $ \Delta_{-n}$ is not empty but  have an empty interior,
    and then $C_{\Delta_{-n}}(X)$ is empty while $C(\Delta_{-n})$ is not. We will also abuse notation concerning
     subsets and write $B\cap C(\Delta_{-n})$ for $\{a \textrm{ restricted to } \Delta_{-n}:a\in B\}$ where
      $B\subset C(X)$.
     \begin{prop}\label{propdeltainvar}   Let $V$ be a closed subset of
$X$ and let $I=C_{X\setminus V}(X)$ be the corresponding ideal. Then for
$n\in \N$ we have
\begin{itemize}
\item[i)] $\al^n(V\cap \Delta_n)\subset V\cap \Delta_{-n}$ iff
$a\in I\cap C(\Delta_{-n}) \Longrightarrow \delta^n(a)\in I \cap  C(\Delta_{n})$,
\item[ii)] $\al^n(V\cap \Delta_n)\supset V\cap \Delta_{-n}$ iff
$\delta^n(a)\in I \Longrightarrow a \in I$, for all
          $a\in C(\Delta_{-n})$.
\end{itemize}
 Hence $V$ is $\al$-invariant ($V\in clos_\al(X)$) if and only if
\begin{equation}\label{n-invariantideal}
 \forall_{n\in\N}\,\,\forall_{a\in C(\Delta_{-n})}
 \,\,a\in I \cap C(\Delta_{-n})\Longleftrightarrow \delta^n(a) \in I\cap C(\Delta_{n}).
\end{equation}
\end{prop}
\begin{Proof} i). Let $\al^n(V\cap\Delta_n)\subset V\cap\Delta_{-n}$ and let $a\in I\cap C(\Delta_{-n})$ be fixed.
Then  for   $x\in V\cap \Delta_{n}$ we have $\al^n(x)\in V\cap\Delta_{-n}$, whence  $\delta^n(a)=a(\al^n(x))=0$ and
  $\delta^n(a)\in I\cap C(\Delta_{n})$.
 \\
  Now suppose $\al^n(V\cap\Delta_n)\nsubseteq V\cap\Delta_{-n}$. Then there exists $x_0\in  V\cap\Delta_n$
   such that $\al^n(x_0)\notin  V\cap\Delta_{-n}$.
  As $\al^n(x_0)\in \Delta_{-n}$ and $V$ is closed, by Urysohn's lemma, there is a function $a_0\in C(X)$
   such that $a_0(\al^n(x_0))=1$ and $a_0(x)=0$ for all $x\in V$. Thus taking $a$ to be the
    restriction of $a_0$ to $\Delta_{-n}$ we obtain $a\in I\cap C(\Delta_{-n})$ but $\delta^n(a)(x_0)=1$,
     whence $\delta^n(a)\notin I\cap C(\Delta_{n})$.\\
ii). Let $\al^n(V\cap\Delta_n)\supset V\cap\Delta_{-n}$ and let
$a\in C(\Delta_{-n})$ be such that $\delta^n(a)\in I\cap C(\Delta_{n})$. Suppose on
the contrary  that $a\notin I\cap C(\Delta_{-n})$.
 Then $a(y_0)\neq 0$ for some $y_0\in V\cap\Delta_{-n}$. Taking $x_0\in V\cap\Delta_n$ such that $y_0=\al^n(x_0)$
 we arrive at the contradiction with $\delta^n(a)\in I\cap C(\Delta_{n})$ because $\delta^n(a)(x_0)=a(y_0)\neq 0$.\\
If $\al^n(V\cap\Delta_n)\nsupseteq V\cap\Delta_{-n}$,  then there exists
$x_0\in V\cap\Delta_{-n}\setminus \al^n(V\cap\Delta_n)$. Similarly as in the proof of item i), using Urysohn's lemma
 we can take $a_0\in C(X)$ such that $a_0(x_0)=1$
  and $a_0|_{\al^n(V\cap\Delta_n)}\equiv 0$. Hence putting  $a=a_0|_{\Delta_{-1}}$
  we have $a\notin I\cap C(\Delta_{-n})$ and $\delta^n(a)\in I\cap C(\Delta_{n})$.\\
In view of  i) and ii),   $\al$-invariance of $V$ is evidently
equivalent to the  condition (\ref{n-invariantideal}).
\end{Proof}
 \begin{defn}
If $I$ is a closed  ideal of $\A$ satisfying
(\ref{n-invariantideal}) then we say that $I$ is \emph{invariant}
under the endomorphism $\delta$, or briefly
\emph{$\delta$-invariant}.
\end{defn}
In virtue of  Proposition \ref{propdeltainvar} it is clear that
$I$ is a $\delta$-invariant ideal iff $I=C_{X\setminus V}(X)$
where $V$ is a closed $\al$-invariant set. Thus, using Theorem
\ref{latticeiso} one can obtain a correspondence between the
invariant ideals of $\A$ and invariant ideals of
$\B$. To this end we denote by $\langle I
\rangle_{_{\B,\,\widetilde{\delta}}}$ the
smallest $\widetilde{\delta}$-invariant ideal of
$\B$ containing $I$.
 \begin{prop}\label{invar-ideals-bijec}
Let $I=C_{X\setminus V}(X)$ be a $\delta$-invariant ideal of $\A$
and let $\V\in  clos_{\tal}(\X)$ be such that $\Phi(\V)=V$. Then
$$
\langle C_{X\setminus V}(X)
\rangle_{_{\B,\,\widetilde{\delta}}} =
C_{\X\setminus \V}(\X),
$$
and the  mapping $I\longmapsto \langle I
\rangle_{_{\B,\,\widetilde{\delta}}}$ establishes
an order preserving bijection between the family of
$\delta$-invariant ideals of $\A$ and $\Tdelta$-invariant ideals
of $\B$. Moreover, the inverse of the mentioned
bijection has the form $\widetilde{I}\longmapsto
\widetilde{I}\cap\A$.
\end{prop}
\begin{Proof} In order to prove the first part of the statement we show that the support
$$S=\bigcup_{f\in\langle C_{X\setminus V}(X)
\rangle_{_{\B,\,\widetilde{\delta}}}}\{\x\in\X:f(\x)\neq
0\}$$
 of the ideal $\langle C_{X\setminus V}(X)
\rangle_{_{\B,\,\widetilde{\delta}}}$ is equal to
$\X\setminus\V$.\\
Let $a\in I=C_{X\setminus V}(X)$. We identify $a$ with $[a,0,...]\in
\B$ and since $x_0\in V$ for any
$\x=(x_0,...)\in\V$, we note that $[a,0,...](\x)=a(x_0)=0$, that
is $a=[a,0,...]\in C_{\X\setminus \V}(\X)$. As $C_{\X\setminus
\V}(\X)$ is $\Tdelta$-invariant, we get $ \langle C_{X\setminus
V}(X) \rangle_{_{\B,\,\widetilde{\delta}}} \subset
C_{\X\setminus \V}(\X)$, whence $S\subset \X\setminus\V$.\\
 Now,
let $\x=(x_0,...,x_k,...)\in \X\setminus \V$. The form of $\V$
(compare Theorem \ref{latticeisomorph}) implies that there exists
$n\in\N$ such that $x_n\notin V$. According to  Urysohn's lemma
there exists $a\in C_{X\setminus V}(X)$ such that $a(x_n)=1$. By
invariance, all the elements $\Tdelta^k(a)$ and $\Tdelta_*^k(a)$
 for $k\in\N$, belong to $\langle I
\rangle_{_{\B,\,\widetilde{\delta}}}$.  In
particular $\Tdelta_*^n(a)=[0,...,a\delta^{n}(1),0,...]\in\langle
I \rangle_{_{\B,\,\widetilde{\delta}}}$ where
$\Tdelta_*(a)(\x)=a(x_n)=1\neq 0$. Thus $\x\in S$ and  
we get $ \X\setminus\V=S $.
\\
In virtue of  Theorem \ref{latticeiso} the relation
$\Phi(\V)=V$ establishes an order preserving bijection between
$clos_{\tal}(\X)$ and $clos_{\al}(X)$ hence the relation $\langle
C_{X\setminus V}(X)
\rangle_{_{\B,\,\widetilde{\delta}}} =
C_{\X\setminus \V}(\X)$ establishes such a bijection too. The
inverse relation $C_{\X\setminus \V}(\X)\cap \A= C_{X\setminus
V}(X)$ is straightforward.
\end{Proof}
 Let us  recall that we
identify  $\B$ with a subalgebra of the covariance
algebra $C^*(\A,\delta)$. Therefore  for any subset $K$ of
$\B$ we denote  by $\langle K \rangle$ an ideal of
$C^*(\A,\delta)$ generated by $K$. The next statement follows from
the preceding proposition and  Theorem 3.5  from
\cite{exel_laca_quigg}.
\begin{thm}\label{idealisomor}
Let $(\A,\delta)$ be a $C^*$-dynamical system such that
 $\al$ has no periodic points.
Then the map
$$
V\longmapsto \langle  C_{X\setminus V}(X) \rangle
$$
is a lattice anti-isomorphism from $clos_\al(X)$  onto the lattice
of ideals in $C^*(\A,\delta)$.  Moreover, for $\V\in
clos_{\tal}(\X)$  such that $\Phi(\V)=V$ the following relations hold
$$
\langle C_{X\setminus V}(X) \rangle=\langle C_{\X\setminus \V}(\X)
\rangle,\qquad \langle C_{X\setminus V}(X) \rangle\cap
\B=C_{\X\setminus \V}(\X),\qquad \langle
C_{X\setminus V}(X) \rangle\cap \A= C_{X\setminus V}(X).
$$
\end{thm}
\begin{Proof} Since $\al$ has no periodic points neither does its reversible extension $\tal$.
The covariance algebra $C^*(\A,\delta)$ is the partial crossed product
$C(\X)\rtimes_{\Tdelta} \Z$ and  $\Z$ is an amenable group. Hence
$(C(\X),\Tdelta)$ has the approximation property, see
\cite{exel_laca_quigg}. Thus in view of  Theorem 3.5 from
\cite{exel_laca_quigg}, the map $\V\longmapsto \langle
C_{\X\setminus \V}(\X) \rangle$ is a lattice anti-isomorphism from
$clos_{\tal}(\X)$  onto the lattice of ideals of $C^*(\A,\delta)$,
and the inverse relation is $\langle C_{\X\setminus \V}(\X)
\rangle\cap C(\X)=C_{\X\setminus
\V}(\X)$.\\
Now, let $\V\in clos_{\tal}(\X)$ and $V\in clos_{\al}(X)$ be such
that $\Phi(\V)=V$. We show  that $\langle C_{X\setminus V}(X)
\rangle=\langle C_{\X\setminus \V}(\X)
\rangle.$\\
On one hand, by Proposition \ref{invar-ideals-bijec} we have
$\langle C_{\X\setminus
\V}(\X) \rangle\cap\A=\langle C_{\X\setminus \V}(\X)\rangle\cap
C(\X)\cap\A= C_{\X\setminus
\V}(\X) \cap\A=C_{X\setminus V}(X)$
and hence $\langle C_{X\setminus V}(X) \rangle\subset \langle
C_{\X\setminus \V}(\X) \rangle$.
On the other hand, $\langle C_{X\setminus V}(X) \rangle\cap C(\X)$
is a $\Tdelta$-invariant ideal of $C(\X)$ containing
$C_{X\setminus V}(X)$, and so it also contains  $\langle
C_{X\setminus V}(X)
\rangle_{_{\B,\,\widetilde{\delta}}} =
C_{\X\setminus \V}(\X)$. Hence $\langle C_{\X\setminus \V}(\X)
\rangle \subset \langle C_{X\setminus V}(X) \rangle$.\\
Concluding,  we have the lattice isomorphism between sets
$V\in clos_\al(X)$ and  $\V\in clos_{\tal}(\X)$, and  the lattice
anti-isomorphism between sets $\V\in clos_{\tal}(\X)$ and ideals
$\langle C_{X\setminus V}(X) \rangle=\langle C_{\X\setminus
\V}(\X) \rangle$, hence $ V\longmapsto \langle  C_{X\setminus
V}(X) \rangle $ is anti-isomorphism.
\end{Proof}
\begin{exm}
In  Example \ref{finitedimension} the only $\al$-invariant sets
are $\{x_1,...,x_{m-1},y_m\}$, $m=1,...,k$, and their sums. The
corresponding ideals are $M_{n_m}$, $m=1,...,k$, and their direct
sums.
\end{exm}
We automatically get a simplicity criteria for the covariance
algebra. We say that $(X,\al)$ forms a cycle, if $X=\{x_0,...,x_{n-1}\}$ is  finite
 and $\al(x_k)=x_{k+1 (mod\,\, n)}$,  $k=0,...n-1$.
\begin{corl}\label{simplicity} Let $\al$ be minimal. If $(X,\al)$ does not form a
cycle then $C^*(\A,\delta)$ is simple.
\end{corl}
\begin{Proof} It suffices to observe that if  $\al$ is minimal then  $\al$
 has no periodic points or $(X,\al)$ forms a cycle.
  Hence we  can  apply Theorem \ref{idealisomor}.
\end{Proof}
\begin{exm}
If $(X,\al)$  does  form a cycle then there are infinitely many
ideals in $C^*(\A,\delta)$.  Indeed if we have $\A=\mathbb{C}^n$ and
$\delta(x_1,...,x_n)=(x_n,x_1,...,x_{n-1})$,  it is known that the
partial crossed product $\mathbb{C}^n\rtimes_\delta\mathbb{Z}_n$
is isomorphic to the algebra $M_n$ of complex matrices $n\times n$ and hence
$C^*(\A,\delta)\hookrightarrow C(S^1)\otimes M_n=C(S^1 ,M_n) $.
\end{exm}
\subsection{The Isomorphism Theorem}
The Isomorphism Theorem simply states that under some conditions epimorphism from
Theorem \ref{UoCa} is in fact an isomorphism. We will prove   here two statements of that kind, Theorems
 \ref{starproperty} and  \ref{isomortheorem},  in the literature however only the latter one is
 named the Isomorphism Theorem. A significant role in proofs of both of these  statements plays
  a certain inequality which ensures the existence of conditional expectation onto the coefficient algebra,
  and which appears in different versions in a number of sources concerning  various crossed products.
  For references  see \cite{lebiedodzij}, \cite{lebied}, \cite{Anton}, \cite{Anton_Lebed},  \cite{Odonovan},
  and for the greatest similarity  with the following Definition \ref{starprop} and Theorem
  \ref{starproperty} see  \cite[Theorem 1.2]{Adji_Laca_Nilsen_Raeburn}.
\begin{defn}\label{starprop}
We  say that a $C^*$-algebra $C^*(\CC,U)$ generated by a  $C^*$-algebra  $\CC$ and an element
$U$ possesses the {\em property\/} $(*)$ if the following
inequality holds
$$
\|\sum_{k=0}^M U^{*k}\pi(a_k^{(0)})U^k\|\leq \|\sum_{n=-N}^N
\Big(\sum_{k=0}^M U^{*k}\pi(a_k^{(n)})U^k\Big)U^n\|\qquad \qquad (*)
$$
for any $a^{(n)}_k\in \CC$ and $M, N\in\N$.
\end{defn}
\begin{thm}\label{starproperty}
Let
$(\pi,U,H)\in\mathrm{CovFaithRep}\,(\A, \delta)$. Then  formula \eqref{kopytko}
establishes   an isomorphism   between the
covariance algebra $C^*(\A,\delta)$  and the $C^*$-algebra
$C^*(\pi(\A),U)$ if and only if $C^*(\pi(\A),U)$ possess the property $(*)$.
\end{thm}
\begin{Proof}
\emph{Necessity}. It suffices to observe that
$C^*(\A,\delta)=C^*(\A,u)$  possess the property $(*)$,  and this
follows
 immediately from the fact that partial crossed products satisfy the appropriate version of this property,
 see  \cite[Remark 2.1]{lebied} and \cite[Proposition 3.5]{mcclanachan}.
\\
\emph{Sufficiency}. By Theorem \ref{coefalgebra2},   $(\pi,U,H)$ extends to the covariant faithful
 representation $(\overline{\pi},U,H)$ of the coefficient $C^*$-algebra $\B$.
 This extended representation satisfies assumptions of \cite[Theorem 3.1]{lebied} and as $\Z$ is amenable
  $(\overline{\pi},U,H)$ give rise to the desired isomorphism, see also \cite[Remark 3.2]{lebied}
\end{Proof}
\begin{corl}\label{suabecorl}
Let $v\in \A$ be a partial isometry such that $uu^*\leq  v^*v,\, vv^*$  where $u$ is the
 universal partial isometry in $C^*(\A,\delta)$.  Then the mapping
$$\Lambda_v(u)=vu, \qquad \Lambda_v(a)=a, \quad a\in \A,
$$
extends to an automorphism of $C^*(\A,\delta)$. In particular,
taking $v=\lambda 1$, $\lambda\in S^1$, we have the  action
$\Lambda$ of the unit circle $S^1$ on $C^*(\A,\delta)$ for which
the fixed points set  is  the coefficient $C^*$-algebra $\B$.
\end{corl}
\begin{Proof} By the above theorem $C^*(\A,\delta)=C^*(\A,u)$  possesses the property $(*)$.
 Clearly, the same is true for $C^*(\A,vu)$.   Since $uu^*\leq  v^*v$ we have   $u=v^*v u=v^*(vu)\in C^*(\A,uv)$,
  whence  $C^*(\A,\delta)=C^*(\A,uv)$, and furthermore $(vu)^*vu=u^*v^*vu=u^*u\in \A$. Since $uu^*\leq vv^*$ we
  have $(vu)a (vu)^*=uau^* vv^*=uau^*$,  that is  the  element $vu$ generates the same endomorphism of $\A$ as $u$,
  and hence applying the preceding theorem  we conclude that $\Lambda_v$ extends to an automorphism of $C^*(\A,\delta)$.
   The rest is straightforward.
\end{Proof}
Now, we are  in position to prove our variant of the  celebrated
Isomorphism Theorem.
\begin{thm}[Isomorphism Theorem]\label{isomortheorem}
Let $(\A, \delta)$  be such that  $\al$ is topologically free.
Then for every $(\pi,U,H)\in\mathrm{CovFaithRep}\,(\A, \delta)$
the algebra $C^*(\pi(\A),U)$ possess property $(*)$.  In other
words, for any two  covariant faithful representations
$(\pi_1,U_1,H_1)$ and $(\pi_2,U_2,H_2)$, the   mapping
$$
U_1\longmapsto U_2,\qquad \pi_1(a)\longmapsto \pi_2(a),\quad a\in \A,
$$
determines an isomorphism of
$C^*(\pi_1(\A),U_1)$ onto $C^*(\pi_2(\A),U_2)$.
\end{thm}
\begin{Proof}
Due to   Theorem \ref{topologicalfreedom}, the  partial
homeomorphism $\tal$ is topologically free and according to 
Theorem \ref{coefalgebra2} representations $\pi_1$ and  $\pi_2$
give rise to covariant representations
$(\overline{\pi}_1,U_1,H_1)$ and $(\overline{\pi}_2,U_2,H_2)$ of
the partial dynamical system $(\B, \Tdelta)$. Thus it
is enough to apply the Theorem 3.6 from \cite{lebied}.
\end{Proof}
\begin{corl}\label{hypotheses}
Let $\A$ act nondegenerately on a Hilbert space $H$,  let
$\delta(\cdot)=U(\cdot)U^*$ where $U\in L(H)$ is a partial
isometry such that $U^*U\in \A$, and let the generated partial
mapping  $\al$ be topologically free.  Then $C^*(\A,U)\cong
C^*(\A,\delta)$.
\end{corl}
The above corollary allow us, in the presence of topological freedom,
 consider only abstract covariance algebras. However, in various concrete specification while using
 the method mentioned
after  Theorem \ref{ideals2},  it may happen that the Isomorphism
Theorem can  be applied to  systems $(\A,\delta)$ such that
$\Delta_{-1}$ is not open and  $\al$ is not topologically free.
\begin{exm}\label{613}
Let $\A$ and $U$ be as in  Example \ref{maxidexm}, then the
associated system forms a cycle and therefore it is not
topologically free.
 However after passing to algebra $\CC=C^*(\A,U^*U)$ we obtain
 the dynamical system $(X\cup\{y\},\al)$ (see Example \ref{maxidexm})\\
\setlength{\unitlength}{.8mm}
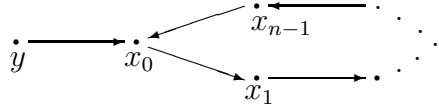
\begin{picture}(80,20)(-55,0)
\put(15,10){\circle*{1}} \put(14,6){$y$} \put(35,10){\circle*{1}}
\put(33,6){$x_0$} \put(55,16){\circle*{1}} \put(54,12){$x_{n-1}$}
\put(75,4){\circle*{1}}
%\put(74,12){$x_{n-2}$}
\put(55,4){\circle*{1}} \put(53,0.5){$x_1$}
\put(53,16){\vector(-3,-1){16}} \put(17,10){\vector(1,0){16}}
\put(37,9){\vector(3,-1){16}} \put(73,16){\vector(-1,0){16}}
\put(57,4){\vector(1,0){16}}
%\put(60,16){.}\put(64,16){.}\put(68,16){.}
\put(74.5,15.5){.}  \put(78,13.7){.} \put(81.5,11.5){.}
\put(78,5){.} \put(81,7){.} \put(84,9.3){.}
\end{picture}\\
which is topologically free, cf. Example \ref{graph}. Hence, due
to the Isomorphism Theorem  $C^*(\A,U)\cong C^*(X\cup\{y\},\al)$. In particular, if  $n=1$ then $C^*(\A,U)$ is the Toeplitz algebra, see Examples 
\ref{exm x_1-x_0}, \ref{exm x_1-x_0-coefalg} and \ref{Toeplitz algebra}.
\end{exm}
\begin{exm}
Consider Hilbert spaces $H_1=L_\mu^2([0,1])$ and
$H_2=L_\mu^2(\R_+)$ where $\mu$ is the Lebesgue measure. We fix  $0<q<1$ and $0<h<\infty$. Let $\A_1\subset L(H_1)$ consists of operators
  of multiplication by
 functions from $C[0,1]$ and let $U_1$  acts according to $(U_1f)(x)=f(q\cdot x)$, $f\in H_1$.
  Similarly, let elements of $\A_2\subset L(H_2)$ act as operators of multiplication by
   functions which are continuous on $\R_+=[0,\infty)$
   and have limit at infinity, and let $U_2$ be the shift operator  $(U_2f)(x)=f( x+h)$, $f\in H_2$.
Then the dynamical systems associated to $C^*$-dynamical systems
$(\A_1,U_1(\cdot)U^{*}_1)$ and $(\A_2,U+2(\cdot)U^{*}_2)$ are
topologically conjugate but the images of the generated mappings are
not open (compare with  Example \ref{q-exm}). Thus  we can not
apply the Theorem \ref{isomortheorem} in the form it is stated.
Nevertheless, endomorphisms of bigger algebras $\CC_1=C^*(\A_1
,U_1^*U_1)$ and $\CC_2=C^*(\A_2,U_2^{*}U_2)$  do generate dynamical
systems
\\
\setlength{\unitlength}{.8mm}
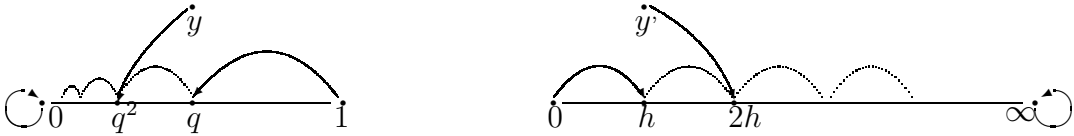
\begin{picture}(180,26)(-17,-2)

%        C[0,1]

\put(5,4){\circle*{1}} \put(6,0){$0$} \put(30,4){\circle*{1}}
\put(29,0){$q$} \put(17.5,4){\circle*{1}} \put(16.5,0){$q^2$}
\put(55,4){\circle*{1}} \put(53.5,0){$1$}
\put(53,4){\line(-1,0){46.5}}

\put(30,20){\circle*{1}} \put(29,16){$y$}

\qbezier(54,5)(42.5,20)(30,5) \put(31,6){\vector(-1,-1){1}}
\qbezier[35](30,5)(23.75,15)(17.5,5)
\qbezier[25](17.5,5)(14.5,11)(11.25,5)
\qbezier[12](11.25,5)(10,8.5)(8.2,5)

\qbezier(29,19.5)(20,12)(17.5,5) \put(18.5,7){\vector(-1,-3){1}}

\put(2,3){\oval(6,6)[b]} \put(2,3){\oval(6,6)[l]}
\put(3.5,5.5){\vector(3,-2){1}}

%      C(R+)

 \put(90,4){\circle*{1}} \put(89,0){$0$}
\put(105,4){\circle*{1}} \put(104,0){$h$} \put(120,4){\circle*{1}}
\put(119,0){$2h$} \put(170,4){\circle*{1}} \put(165,1){$\infty$}

\put(168,4){\line(-1,0){76.5}}

\put(105,20){\circle*{1}} \put(103.5,16){$y^,$}
\qbezier(106,20)(115.5,15)(119.5,6)\put(119.1,7.4){\vector(1,-3){1}}

\qbezier(90,5)(97.5,15)(104.5,5.4)
\put(104.2,6.7){\vector(1,-2){1}}
\qbezier[40](105,5)(112.5,15)(119.5,5)
\qbezier[33](120,5)(127.5,15)(134.5,5)
\qbezier[25](136,5)(142.5,15)(149.5,5)

\put(173,3){\oval(6,6)[r]} \put(173,3){\oval(6,6)[bl]}
\put(172,5.5){\vector(-2,-1){1}}
\end{picture}\\
satisfying the assumptions of the Isomorphism Theorem. These
dynamical systems are topologically conjugate by a piecewise
linear mapping $\phi$  which maps
 $nh$ into $q^n$, $n\in \N\cup\{\infty\}$, and $y'$ into $y$, that is
$$
\phi(x)=q^n\big(\frac{q-1}{h}\,x+ 1-n(q-1)\big),\qquad \mathrm{for}\,\, x\in [nh,nh+1),\qquad \mathrm{and}\,\,\,\,\,\,\phi(\infty)=0,\,\,\, \phi(y')=y.
$$
Therefore, by the Isomorphism Theorem, the mapping $\A_1\ni a\mapsto
a\circ\phi\in \A_2$, and $U_1\mapsto U_2$, establishes the expected
isomorphism; $ C^*(\A_1,U_1)\cong C^*(\A_2,U_2).$
\end{exm}

Lastly, we would like to present an example which shows how the  results achieved in this paper clarify  the situation mentioned in the example from which we have started  the introduction.
\begin{exm}[Solenoid]\label{solenoid2}
Let $H=L_\mu^2(\R)$ where $\mu$ is the Lebesgue measure on  $\R$, and let $\A\subset L(H)$ consists of the operators of multiplication by periodic continuous functions  with period $1$, that is  $\A\cong C(S^1).$ 
Set the unitary operator  $U\in L(H)$ by the formula 
$$(Uf)(x)=\sqrt{2}\,f(2x).
%\qquad  (U^*f)(x)= \frac{1}{\sqrt{2}}\,f\Big(\frac{x}{2}\Big).$$   
$$
Then for each $a(x)\in \A$, $UaU^*$ is the  operator of mulitplication by the periodic function  $a(2x)$ with period $\frac{1}{2}$, and
$U^*aU$ is the operator of multiplication by  the periodic function  $\displaystyle{a\Big(\frac{x}{2}\Big)}$ with period $2$.  Hence  
$$
U\A U^*\subset \A\,\,\,\, \textrm{ and  }\,\,\,\, U^* \A U\nsubseteq \A.
$$
The endomorphism $U(\cdot)U^*$ generate on the spectrum of $\A$ the mapping $\al$ given by $\al(z)=z^2$ for $z\in S^1$, and the spectrum  of the algebra $\B$ generated by  $\bigcup_{n\in\N} U^{*n}\A U^{n}$ is the solenoid $\Sol$: $\B\cong C(\Sol)$, cf. Example \ref{solenoid}.  Furhetmore  $\al$ is topologically free and therefore we have
$$
C^*(\A,U) \cong C^*(S^1,\al)=C(\Sol) \rtimes_F \Z
$$
where in the right hand side stands the standard crossed product of $\B=C(\Sol)$ by the {\em automorphism} induced by  the solenoid map $F$, see  Example \ref{solenoid}.

\end{exm}
\section*{Summary}
\addcontentsline{toc}{section}{Summary}

In this paper we introduced crossed product-like realization of
the universal algebra associated to 'almost' arbitrary commutative
$C^*$-dynamical system $(\A,\delta)$. This new realization
generalizes the known constructions for $C^*$-dynamical systems
where dynamics is implemented by an automorphism or a
monomorphism.
\\
  The main gain of this  is that we  are able to describe  important characteristics of the
  investigated object  in terms of the underlying
topological (partial) dynamical system $(X,\al)$, the tool which
until now was  used successfully  only in the case of a (partial)
automorphism.
\\
Namely, we have described the ideal structure of covariance algebra by closed invariant subsets of $X$,
 in particular simplicity criteria is obtained.
Moreover we have generalized  the topological freedom, the
condition  under which all the covariant  faithful representations
of $(\A,\delta)$ are algebraically equivalent, see the Isomorphism Theorem. For applications this is probably the most important result of the  paper. 
\\
The important novelty in our approach is that the  construction of
covariance algebra  here consists of two independent steps. The
advantage of this is that one may analyse covariance algebra on
two levels. First, one may study the relationship between initial
$C^*$-dynamical system and the one generated on its coefficient
$C^*$-algebra, and then one may apply known statements and methods
as the latter system is more accessible (generated mapping on the
spectrum of coefficient algebra is bijective).
\\
We indicate that recently (see \cite{exel2}) a notion of crossed-product
of  a $C^*$-algebra by  an endomorphism (or even partial endomorphism, see \cite{exel_royer})
 has been introduced, a construction which depends also on the choice of transfer operator.
  This construction is especially  well adapted to deal with morphisms which generate local homeomorphisms.
  In particular it was used to investigate Cuntz-Krieger algebras,
  cf. \cite{exel2}, \cite{exel_royer}, \cite{exel_laca}.
   However it seems that in case $\al$ is not injective  there does not
   exist a transfer operator such that the aforementioned crossed-product
   is isomorphic to covariance algebra considered here, and in case $\al$ is injective
    the transfer operator is trivial, that is, it is $\al^{-1}$ and thus it does not add anything new to the system.

\bigskip

\centerline{\textsc{Acknowledgements}}
\bigskip

\noindent The author  wishes to express his thanks  to  A. V. Lebedev  for suggesting
  the problem and many stimulating conversations,  to  A. K. Kwa\'sniewski for his active
   interest in the preparation of this paper, and also to D. Royer for pointing an error in
   an earlier version of Theorem 2.2.

\end{document}